\documentclass[12pt,titlepage]{amsart}

\usepackage{amsmath}
\usepackage{amsthm}
\usepackage{amssymb}
\usepackage{mathabx}
\usepackage{mathrsfs} 
\usepackage[all]{xy} 

\usepackage{thmtools}    
\usepackage{thm-restate} 

\usepackage{graphicx}
\graphicspath{{figures/}} 

\usepackage{quickmath} 

\usepackage{cancel} 
\usepackage{enumitem}
\usepackage{hyperref} 
\usepackage[usenames,dvipsnames]{color} 
\usepackage[marginpar]{todo}

\usepackage{tikz}
\usetikzlibrary{arrows}
\usetikzlibrary{calc}

\hypersetup{colorlinks=true, linkcolor=NavyBlue, urlcolor=WildStrawberry, citecolor=blue}

\title[Recognizing RACGs Using Involutions]{Recognizing Right-Angled Coxeter Groups Using Involutions}

\author[Cunningham, Eisenberg, Piggott, Ruane]{%
Charles Cunningham \\ Tufts University \\ {~} \\
Andy Eisenberg \\ Tufts University \\ {~} \\
Adam Piggott \\ Bucknell University \\ {~} \\
Kim Ruane \\ Tufts University \\ (corresponding author: \texttt{kim.ruane@tufts.edu}){~} \\ {~} \\ 
Keywords: Coxeter group, involutions, graph theory, automorphisms}

\let\oldtocsection=\tocsection
\let\oldtocsubsection=\tocsubsection
\renewcommand{\tocsection}[2]{\hspace{0em}\oldtocsection{#1}{#2}}
\renewcommand{\tocsubsection}[2]{\hspace{1.1em}\oldtocsubsection{#1}{#2}}

\newcommand{\Inv}[1][]{\Delta_{#1}}

\newcommand{\Ab}[1]{#1^{\textrm{ab}}}

\begin{document}

\maketitle

\begin{abstract}
We consider the question of determining whether or not a given group (especially one generated by involutions) is a right-angled Coxeter group.  We describe a group invariant, the \emph{involution graph}, and we characterize the involution graphs of right-angled Coxeter groups.  We use this characterization to describe a process for constructing candidate right-angled Coxeter presentations for a given group or proving that one cannot exist.  We apply this process to a number of examples.  Our new results imply several known results as corollaries.  In particular, we provide an elementary proof of rigidity of the defining graph for a right-angled Coxeter group, and we recover an existing result stating that if $\Gamma$ satisfies a particular graph condition (called \emph{no SILs}), then $\Aut^0(W_{\Gamma})$ is a right-angled Coxeter group.
\end{abstract}

\tableofcontents

\section{Introduction}

Given a finite simple graph $\Gamma$, the \emph{right-angled Coxeter group defined by $\Gamma$} is the group $W = W_{\Gamma}$ generated by the vertices of $\Gamma$.  The relations of $W_{\Gamma}$ declare that the generators all have order 2, and adjacent vertices commute with each other. Right-angled Coxeter groups (commonly abbreviated RACG) have a rich combinatorial and geometric history \cite{Davis}. The particular presentation specified by $\Gamma$ is called a \emph{right-angled Coxeter system}.  When encountering a group generated by involutions, a natural question is to ask whether or not this group might be a right-angled Coxeter group, and if so, how to identify the preferred presentation.

The main objective of this paper is the development of a recognition procedure that successfully answers this question for certain families of groups.  Although the procedure may be applied more generally, our applications focus primarily on two classes of examples.  Given a right-angled Coxeter group $W_{\Gamma}$, we consider
\enum
\item extensions of $W_{\Gamma}$ by subgroups of $\Out^0(W_{\Gamma})$, and
\item subgroups of $W_{\Gamma}$ generated by chosen sets of involutions.
\eenum
(Recall that $\Aut^0(W_{\Gamma})$ consists of the automorphisms of $W_{\Gamma}$ which map each generator to a conjugate of itself, and $\Out^0(W_{\Gamma})$ is the quotient $\Aut^0(W_{\Gamma})/\Inn(W_{\Gamma})$.)  In each of these cases, we give examples of groups which are right-angled Coxeter and examples which are not.    For those cases which are right-angled Coxeter, our procedure produces the preferred presentations.   We show:

\begin{restatable*}{Theorem}{ExtensionByPCs}
\label{thm:mainextensiontheorem}
Suppose $\chi_1, \dots, \chi_k$ are pairwise commuting partial conjugations of the right-angled Coxeter group $W_{\Gamma}$ such that whenever $\chi_i$ and $\chi_j$ have the same acting letter, their domains don't intersect.  Then $G = W \rtimes \langle \chi_1, \dots, \chi_k \rangle$ is a right-angled Coxeter group.  Further, writing $S_i \subseteq \{\chi_1, \dots, \chi_k\}$ for the set comprising those partial conjugations with acting letter $a_i$, we have that
\[
\left\{a_1 \prod_{\chi_i \in S_1} \chi_i, \dots, a_n \prod_{\chi_i \in S_n} \chi_i\right\} \cup \left\{\chi_1, \dots, \chi_k\right\}
\]
is a Coxeter generating set for $G$.
\end{restatable*}


If a group $G$ has only 2-torsion, and $G$ is not a right-angled Coxeter group, then $G$ is not a Coxeter group.  So our procedure may in fact enable one to show that a given group is not a Coxeter group.  The first author has used some of the methods described here to show that $\Out^0(W_n)$ for $n \geq 4$ is not a Coxeter group \cite{Charlie}.  ($W_n$ is the \emph{universal Coxeter group} whose defining graph has $n$ vertices and no edges.)  

Given a group $G$, the \emph{involution graph $\Inv[G]$ of $G$} is the group invariant defined as follows: the vertices in $\Inv[G]$ correspond to the conjugacy classes of involutions in $G$; vertices are adjacent when there exist commuting representatives of the corresponding conjugacy classes.
In general, this invariant is unwieldy.  It may be infinite, and even when it's finite, it may be impossible to construct.  Nevertheless, for certain classes of groups the invariant promises insights. Like any invariant, it can allow us to distinguish between groups.  It also carries information on the automorphism group of $G$.
Since an automorphism must permute conjugacy classes of involutions and must preserve commuting relations, $\Aut(G)$ acts naturally on $\Inv[G]$.  The kernel of this action is therefore a natural normal subgroup of $\mathrm{Aut}(G)$, and has finite index in $\Aut(G)$ when $\Inv[G]$ is finite. 

The involution graph for a right-angled Coxeter group $W_{\Gamma}$ is easily constructed directly from $\Gamma$: the vertices in $\Inv[W]$ correspond to cliques in $\Gamma$; vertices are adjacent when the union of the corresponding cliques is also a clique.  When constructed in this manner, we denote the graph $\Gamma_K$ and call it the \emph{clique graph for $\Gamma$}.  Tits \cite{Tits} proved that the kernel of the action $\Aut(W) \circlearrowright \Inv[W]$ has a natural complement, which is therefore a finite subgroup of $\Aut(\Inv[W])$.  Thus the involution graphs of right-angled Coxeter groups are significantly more tractable than the involution graphs of arbitrary groups, and may be more convenient for certain purposes than the defining graph $\Gamma$.  Aaron Meyers, in his undergraduate thesis under the supervision of the third author, began to explore some properties of clique graphs and how to recover their base graphs.  (As this work is unpublished, new proofs are given in the following sections.)

The reader may compare our use of the clique graph and involution graph to the use of the \emph{clique graph}, \emph{extension graph}, and \emph{commutation graph} in \cite{KK} in the context of right-angled Artin groups.  Our use of the term and notation for the clique graph comes from \cite{KK}.  In addition, Kim and Koberda define the \emph{extension graph $\Gamma^e$ of $\Gamma$} and the \emph{commutation graph} of a subset $S \subset A(\Gamma)$ of elements of the right-angled Artin group.  The vertices of $\Gamma^e$ are the words in the right-angled Artin group $A(\Gamma)$ which are conjugate to a vertex of $\Gamma$, and two such vertices are connected by an edge if they commute with one another.  More generally, the commutation graph of $S$ has vertices given by the elements of $S$, and two of these are connected by an edge if they commute with each other.

It is straightforward to define the extension and commutation graphs in the context of right-angled Coxeter groups.  Note that the vertices of $\Gamma^e$ are the individual group elements, not conjugacy classes, so that $\Gamma^e$ is infinite whereas $\Inv[W_{\Gamma}]$ is finite.  Moreover, $\Gamma^e$ does not contain words that are only conjugate to a product of pairwise commuting generators, so it is not the case that the $\Inv[W_{\Gamma}]$ is a quotient graph of $\Gamma^e$.  Theorem 1.3 in \cite{KK} states that, given graphs $\Lambda$ and $\Gamma$, if $\Lambda$ is contained in $\Gamma^e$, then $A(\Lambda) \leq A(\Gamma)$.  The analogous statement about right-angled Coxeter groups is certainly false, and a counterexample is provided by
\[
D_{\infty} = W_2 = \langle a, b \mid a^2 = b^2 = 1 \rangle.
\]
The defining graph $\Gamma$ consists of two vertices with no edges.  The extension graph $\Gamma^e$ has countably many vertices and no edges, but $D_{\infty}$ cannot contain subgroups which are free products of more than two copies of $\Z/2\Z$.  If we replace the extension graph with the involution graph $\Inv[W_{\Gamma}]$ in Theorem 1.3 in \cite{KK}, the claim would still be false.  $\Inv[W_{\Gamma}]$ contains cliques which are larger than any clique in $\Gamma$.  

Finally, we note that the involution graph $\Inv[G]$ of a group which is not a right-angled Coxeter group may not be a commutation graph on any subset $\{g_1, \dotsc, g_n\}$ of elements.  \emph{A priori}, it could be the case that there is no single collection of elements, one from each conjugacy class, which simultaneously exhibit all commuting and non-commuting relationships dictated by the involution graph.  (In the case of a right-angled Coxeter group $W_{\Gamma}$, $\Inv[W_{\Gamma}]$ is the commutation graph on the set of products of pairwise commuting generators.)  It may be that the techniques of \cite{KK} could be adapted to the case of right-angled Coxeter groups, but as the current paper focuses on the recognition problem, the authors have not considered questions of embedability.

In Section \ref{sec:summary}, we summarize our recognition procedure which attempts to construct right-angled Coxeter presentations for a given group.  This procedure relies on many facts about clique graphs and involution graphs which, for clarity of exposition, are only stated in that section.  Detailed proofs have been relegated to Section \ref{sec:details} at the end of the paper.  Section \ref{sec:summary} contains all necessary definitions and results to understand the applications in Section \ref{sec:applications}.

In Section \ref{sec:applications}, we apply our procedure to several first examples of potential right-angled Coxeter groups.  Section \ref{sec:positiveresults} collects examples of families of groups which are right-angled Coxeter.  $\Gamma$ is said to contain a \emph{separating intersection of links (SIL)} if, for some pair of vertices $v$ and $w$ with $d(v, w) \geq 2$, there is a connected component of $\Gamma \setminus (\Lk(v) \cap \Lk(w))$ which contains neither $v$ nor $w$.  Otherwise, we say \emph{$\Gamma$ contains no SILs}.  Section \ref{sec:positiveresults} also gives a new, shortened proof of a prior result \cite[Theorem 3.6]{CRSV}: that $\Aut^0(W_{\Gamma})$ is right-angled Coxeter if $\Gamma$ contains no SILs.  Section \ref{sec:negativeresults} shows several examples of groups which we prove cannot be right-angled Coxeter.  This includes, in particular, an iterated extension
\[
\overbrace{\underbrace{\left(W_{\Gamma}\rtimes \Z/2\Z\right)}_{H} \rtimes \Z/2\Z}^{G}
\]
in which $H$ is not right-angled Coxeter, but $G$ is.  We also note that $\Aut^0(W_3)$ is not right-angled Coxeter, answering a motivating question for the authors.

Section \ref{sec:decompositionresults} states some results that essentially identify features of a given graph $\Lambda$ which indicate that $W_{\Lambda}$ has a semi-direct product decomposition $W_{\Lambda} = W_{\Gamma} \rtimes H$, where $H \leq \Out^0(W_{\Gamma})$.  The results of this section follow from those in Section \ref{sec:positiveresults} quite easily, and the semi-direct product decompositions are certainly not unique.

Section \ref{sec:details} presents detailed proofs about many facts stated without proof in Section \ref{sec:summary}.  In this section, we present a characterization of those finite graphs which arise as clique graphs (i.e., a characterization of those graphs which arise as the involution graphs of right-angled Coxeter groups).  We present a collapsing procedure to recover $\Gamma$ from $\Gamma_K$, and we establish the correctness of our recognition procedure for constructing right-angled Coxeter presentations.

Finally, in Section \ref{sec:future} we give many follow-up questions which may be approachable using our recognition procedure.  These include the question of characterizing those subgroups $H \leq \Out^0(W_{\Gamma})$ such that $W_{\Gamma} \rtimes H$ is again right-angled Coxeter, and determining when the involution graph of a subgroup $H \leq G$ can be calculated easily from the involution graph of $G$.

\section{A Summary of the Recognition Algorithm}
\label{sec:summary}

In this section, we present the definitions and basic properties of the clique graph, star poset, and involution graph constructions.  We state one of our main theorems characterizing those finite graphs which arise as clique graphs, and we describe a procedure which recovers a graph $\Gamma$ from its clique graph $\Gamma_K$.  Finally, we prove several algebraic results about right-angled Coxeter groups which allow us to modify this procedure to seek right-angled Coxeter presentations of a given group.  Many of the proofs of this section are elementary or non-geometric in nature, so they have been pushed to Section \ref{sec:details} at the end of the paper, where the interested reader will find all of the details.  In this section, we present only the definitions and statements of results necessary to understand the applications in Section \ref{sec:applications}.

A finite simple graph $\Gamma = (V, E)$ is an ordered pair of finite sets.  We require that $V$, the set of \emph{vertices}, is nonempty and $E$, the set of \emph{edges}, consists of 2-element subsets of $V$.  We say $a, b\in V$ are \emph{adjacent} if $\{a, b\} \in E$.  All graphs we consider in this paper will be undirected and have finitely many vertices, no loops, and no parallel edges.  We will use the notations
\begin{align*}
\Lk(v) &= \{w\in V \mid \{v, w\} \in E\} \\
\St(v) &= \Lk(v) \cup \{v\}
\end{align*}
for the \emph{link of $v$} and the \emph{star of $v$}, respectively.

\defi
Let $\Gamma$ be a graph.  A \emph{clique in $\Gamma$} is a nonempty subset of pairwise adjacent vertices.  The \emph{clique graph of $\Gamma$} is the graph $\Gamma_K = (V_K, E_K)$ whose vertices correspond to the cliques of $\Gamma$.  Two vertices of $\Gamma_K$ are adjacent if the union of the corresponding cliques in $\Gamma$ is also a clique.
\edefi

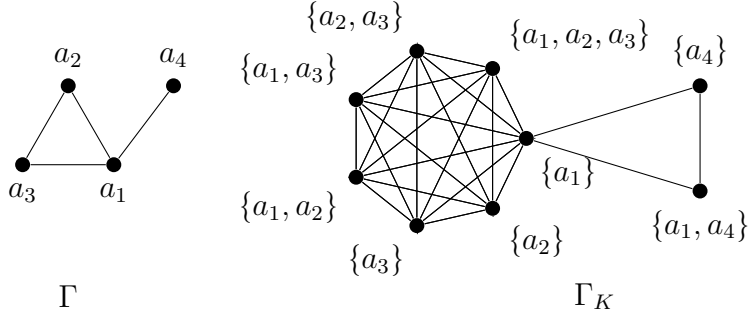
\begin{figure}[ht]
\begin{tikzpicture}[scale=0.7]
\node (a1) [circle,fill,inner sep=2pt,label=-90:$a_1$] at (-30:1) {};
\node (a2) [circle,fill,inner sep=2pt,label=90:$a_2$] at (90:1) {};
\node (a3) [circle,fill,inner sep=2pt,label=-90:$a_3$] at (210:1) {};
\node (a4) [circle,fill,inner sep=2pt,label=90:$a_4$] at (2,1) {};

\draw (a1) to (a2) to (a3) to (a1) to (a4);

\node at (0, -3) {$\Gamma$};

\begin{scope}[xshift=7cm]
\node (b1) [circle,fill,inner sep=2pt,label={1*360/7}:{$\{a_1, a_2, a_3\}$}] at ({1*360/7}:1.7) {};
\node (b2) [circle,fill,inner sep=2pt,label={2*360/7}:{$\{a_2, a_3\}$}] at ({2*360/7}:1.7) {};
\node (b3) [circle,fill,inner sep=2pt,label={3*360/7}:{$\{a_1, a_3\}$}] at ({3*360/7}:1.7) {};
\node (b4) [circle,fill,inner sep=2pt,label={4*360/7}:{$\{a_1, a_2\}$}] at ({4*360/7}:1.7) {};
\node (b5) [circle,fill,inner sep=2pt,label={5*360/7}:{$\{a_3\}$}] at ({5*360/7}:1.7) {};
\node (b6) [circle,fill,inner sep=2pt,label={6*360/7}:{$\{a_2\}$}] at ({6*360/7}:1.7) {};
\node (b7) [circle,fill,inner sep=2pt,label={-70}:{$\{a_1\}$}] at ({7*360/7}:1.7) {};

\node (b8) [circle,fill,inner sep=2pt,label={90}:{$\{a_4\}$}] at (5,1) {};
\node (b9) [circle,fill,inner sep=2pt,label={-90}:{$\{a_1, a_4\}$}] at (5,-1) {};

\foreach \x in {1, 2, 3, 4, 5, 6, 7}
    {
    \foreach \y in {1, 2, 3, 4, 5, 6, 7}
        {
            \draw (b\x) to (b\y);
        }
    }

\draw (b7) to (b8) to (b9) to (b7);

\end{scope}
\node at (10, -3) {$\Gamma_K$};

\end{tikzpicture}
\caption{An example of a graph $\Gamma$ (left) and its corresponding clique graph $\Gamma_K$ (right).}
\label{fig:cliqueexample}
\end{figure}

The relation $v \sim w$ when $\St(v) = \St(w)$ is an equivalence relation on $V(\Gamma)$.  Write $[v]$ for the equivalence class of $v$.  Declaring $[v] \leq [w]$ if $\St(v) \subseteq \St(w)$ defines a partial ordering, and we write $\mathcal{P}(\Gamma)$ for the poset of star-equivalence classes of vertices in $\Gamma$.

Throughout this paper, we will write $\Gamma_1, \Gamma_2, \dotsc, \Gamma_r$ for the maximal cliques of $\Gamma$.  If $I \subset \{1, 2, \dotsc, r\}$, then
\[
\Gamma_I = \bigcap_{i\in I} \Gamma_i
\]
is the corresponding intersection of maximal cliques.

\defi
A vertex $v\in \Gamma$ is called \emph{minimal} if it is contained in a unique maximal clique.  Given $J \subset \{1, 2, \dotsc, r\}$, we say $v$ is \emph{$J$-minimal} if there is no $J' \supset J$ such that $\Gamma_{J'} \subsetneq \Gamma_J$ and $v\in \Gamma_{J'}$.
\edefi

\begin{restatable}{Theorem}{cliquecharacterization}
\label{thm:cliquecharacterization}
Let $\Gamma'$ be a graph.  Then there exists a graph $\Gamma$ such that $\Gamma' = \Gamma_K$ if and only if the following three conditions are satisfied:
\enum
\item (Maximal Clique Condition) For all $I$, there exists some $k_I$ such that
\[
|\Gamma_I'| = 2^{k_I} - 1.
\]
\item (Minimal Vertex Condition) Each nonempty intersection $\Gamma_J'$ contains some $J$-minimal vertex $v_J$.
\item (Inclusion-Exclusion Condition) For each $J$, 
\[
\sum_{I \supsetneq J} (-1)^{|I \setminus J| + 1} k_I \leq k_J.
\]
\eenum
\end{restatable}

Moreover, if $\Gamma'$ is a clique graph, then the graph $\Gamma$ such that $\Gamma' = \Gamma_K$ is unique.  The following procedure, which we call the \emph{collapsing procedure}, recovers $\Gamma$ from $\Gamma'$.  We may write $\Gamma = C(\Gamma')$.

\begin{restatable}{Theorem}{collapsing}
\label{thm:collapsing}
Let $\Gamma'$ be a graph which satisfies the Maximal Clique, Minimal Vertex, and Inclusion-Exclusion Conditions.  Then there is a unique (up to isomorphism) graph $\Gamma$ such that $\Gamma'$ is isomorphic to $\Gamma_K$.  Moreover, the following \emph{collapsing procedure} produces the graph $\Gamma$ if it exists.
\enum[labelindent=\parindent,label=\arabic*:]
\item Initially, let $V = \{\}$.
\item \label{step:choosing} Let $[w] \in \mathcal{P}(\Gamma')$ be a class such that every class $[v]$ with $[w] < [v]$ has already been considered.  Write
\[
S_w = \bigcup_{[v] \geq [w]} [v].
\]
Then there is some $k$ such that $|S_w| = 2^k - 1$.  Let $k'$ be the number of vertices of $S_w$ which are already contained in $V$.  Choose $k - k'$ vertices of $[w]$ to add to the vertex set $V$.
\item Repeat the previous step until all classes of $\mathcal{P}(\Gamma')$ have been considered.
\item Return the graph $C(\Gamma')$ which is the induced subgraph of $\Gamma'$ on the vertex set $V$.
\eenum
\end{restatable}

We remark that the set $S_w$ forms a clique in $\Gamma'$ which is an intersection of maximal cliques, so its size has the desired form by the Maximal Clique Condition.  The details can be found in Section \ref{sec:cliquegraph}.

\defi
Let $G$ be a (finitely generated) group.  The \emph{involution graph of $G$}, denoted, $\Inv[G]$, is a graph defined as follows.  The vertices are the conjugacy classes of involutions in $G$.  Two vertices $[x]$ and $[y]$ are connected by an edge if there exist representatives $gxg^{-1}$ and $hyh^{-1}$ which commute with each other.
\edefi

We make a few remarks.  The particular conjugates which witness commutativity are chosen for each edge individually.  A system of representatives of each conjugacy class which act as witnesses for every edge simultaneously is called a \emph{full system of representatives}.  Such a system need not exist in general, but a right-angled Coxeter group will always have a full system of representatives.

We have also said earlier that all graphs we consider do not have loops, although the involution graph as defined here may contain a loop if an involution commutes with a conjugate of itself.  This may happen in general, but it will never happen in a right-angled Coxeter group.  So, if the involution graph of a given $G$ contains a loop, we may immediately conclude that $G$ is not a right-angled Coxeter group.

\lem
Let $\Gamma$ be a graph.  Then $\Inv[W_{\Gamma}] = \Gamma_K$.
\elem
\pf
It is a well-known fact about right-angled Coxeter groups that the only nontrivial torsion elements have order 2, and that any involution is conjugate to some product of pairwise commuting generators.  The set of products of pairwise commuting generators forms a full system of representatives for the involution graph (this follows essentially from the deletion condition), and two such products commute if and only if all the generators involved in each product pairwise commute (i.e., if the collection of all these generators forms a clique in $\Gamma$).
\epf

We recover the rigidity of right-angled Coxeter groups as an immediate consequence.  This was originally proven in \cite{Green} (for a more general class of groups), and many other proofs have been presented for different classes of groups containing right-angled Coxeter groups as a subclass (see, for examples, \cite{Droms,Laurence,Radcliffe}).

\cor
The defining graph of a right-angled Coxeter group $W_{\Gamma}$ is unique up to isomorphism.
\ecor
\pf
The involution graph is an algebraic invariant (it does not depend on the chosen right-angled Coxeter presentation).  By the previous lemma, the involution graph $\Inv[W_{\Gamma}]$ is a clique graph, and by Theorem \ref{thm:cliquecharacterization} the collapsed graph $C(\Inv[W_{\Gamma}])$ is unique (up to isomorphism).
\epf

At this point, we can essentially describe our recognition procedure for seeking a right-angled Coxeter presentation for a given group $G$.  First, we form the involution graph $\Inv[G]$.  If this is not a clique graph, then $G$ is not a right-angled Coxeter group.  If it is, then we must find a full system of representatives for the vertices.  If such a system does not exist, then $G$ is not a right-angled Coxeter group.  If we find a full system of representatives, then the collapsing procedure will produce a labeled graph $\Gamma = C(\Inv[G])$, which gives a map $W_{\Gamma} \to G$ by sending the generators of $W_{\Gamma}$ to the labels of the corresponding vertices.  If we can show the candidate map is an isomorphism, then $G$ is a right-angled Coxeter group, and the labels of $\Gamma$ form a right-angled Coxeter generating set.  (On the other hand, if the candidate map is not an isomorphism, we cannot conclude that $G$ is not a right-angled Coxeter group.  We may have simply chosen the wrong full system of representatives for $\Inv[G]$.)

We must address one subtlety in this procedure.  In Theorem \ref{thm:collapsing}, we chose vertices from $[w]$ to add to the vertex set $V$ arbitrarily.  It only mattered that we had the right number of vertices from each intersection of maximal cliques.  In the algebraic setting, this is not sufficient, as the following simple example shows.

\ex
Let $\Gamma$ be a triangle with vertices $a, b, c$.  Then $\Gamma_K = \Inv[W_{\Gamma}]$ is a clique of size 7 with the labels $a, b, c, ab, ac, bc, abc$.  In the star poset $\mathcal{P}(\Gamma_K)$, all vertices are equivalent, so there is only one $[w]$ to consider.  The collapsing procedure says to choose 3 vertices from this class at random.  If we choose, for example, the vertices $a, b, c$, then the collapsing procedure recovers $\Gamma$.  If we choose $a, ab, abc$, then we find a new right-angled Coxeter presentation for $W_{\Gamma}$.  However, if we pick $a, b, ab$, then we don't get a right-angled Coxeter presentation (because there is an additional relation between these vertices).
\eex

Essentially, at this step in the collapsing procedure we are choosing which vertices of the involution graph represent generators and which represent products of generators.  There are (generally) many different ways that we can make this choice, but we have to make use of some algebraic information to avoid choosing products as if they were generators.  The following results are certainly of independent interest, but we will, in particular, use them to make intelligent choices during the collapsing procedure.

Since we wish to avoid choosing vertices whose labels have a nontrivial product relation, it would certainly help if we could solve the word problem in $G$.  However, depending on how $G$ is presented, such a solution may or may not be evident (if it even exists).  For this reason, we will pass to the abelianization $\Ab{G}$, in which there is a solution to the word problem.  If $G$ is a right-angled Coxeter group, then $\Ab{G} \cong (\Z/2\Z)^n$, and a product relation among involutions in $G$ must also occur in $\Ab{G}$.

From this point forward, for $g\in G$, we will write $\overline{g}$ for the image of $g$ in the abelianization.  An important fact about right-angled Coxeter groups is that the abelianization is injective on conjugacy classes of involutions.

\prop 
\label{prop:abel}
Let $W_{\Gamma}$ be a right-angled Coxeter group.  Let $x,y \in W_{\Gamma}$ such that $x^2 = y^2 = 1$. Then $\overline{x} = \overline{y}$ in $\Ab{W_{\Gamma}}$ if and only if $x$ and $y$ are conjugate in $W_{\Gamma}$.  
\eprop
\pf
The ``if'' direction is trivial.  Now, suppose that $x$ and $y$ are not conjugate in $W_{\Gamma}$. Since $x, y$ are involutions, there are pairwise commuting generators $a_1, a_2, \dotsc, a_k$, pairwise commuting generators $b_1, b_2, \dotsc, b_{\ell}$, and words $g, h$ such that
\begin{align*}
x &= g a_1 a_2 \dotsm a_k g^{-1} \\
y &= h b_1 b_2 \dotsm b_{\ell} h^{-1}.
\end{align*}
Without loss of generality, since $x$ and $y$ are not conjugate, there is a $b_j$ that does not appear among the $a_i$. But since it is a generator, there is a $\Z/2\Z$ direct factor in $\Ab{W_{\Gamma}}$ corresponding to that $\overline{b_j}$. Therefore, $\overline{y}$ will have a 1 in this factor and $\overline{x}$ will have a 0. Thus, $\overline{x} \neq \overline{y}$ in $\Ab{W_{\Gamma}}$.
\epf

\cor
\label{cor:inj}
For a right-angled Coxeter group $W_{\Gamma}$, if $H$ is a subgroup generated by distinct, commuting involutions, then $H \cong \Ab{H}$ injects into $\Ab{W_{\Gamma}}$. 
\ecor
\pf
$H$ is a finite subgroup of $W_{\Gamma}$ and so is conjugate to a special subgroup $H'$. Each element of $H'$ is a distinct product of commuting generators from $W_{\Gamma}$ and so each gets sent to a distinct element of $\Ab{W_{\Gamma}}$. Thus, no two elements of $H'$ can be conjugate in $W_{\Gamma}$ and so neither can any two elements of $H$. By Proposition \ref{prop:abel}, $H$ injects into $\Ab{W_{\Gamma}}$. 
\epf

\begin{restatable}{Proposition}{propchoice}
\label{prop:choose}
If $W_{\Gamma}$ is a right-angled Coxeter group, then in step 2 of the collapsing procedure in Theorem \ref{thm:collapsing}, we can choose the $k-k'$ involutions of $W_{\Gamma}$ so that the chosen elements do not exhibit a non-trivial product relation.
\end{restatable}
\pf
See Section \ref{sec:involutiongraph}.
\epf

This proposition makes use of the available algebraic information to amend our collapsing procedure and avoid nontrivial product relations.  We can make further use of the available algebraic information to improve upon the procedure.  In general, we have no particular method (or hope of finding a method) to construct $\Inv[G]$ for an arbitrary $G$.  Each of the following steps seem to be generally insurmountable:
\enum
\item Identify all involutions in $G$.
\item Separate all involutions into their conjugacy classes.
\item Determine the presence or lack of each edge in $\Inv[G]$ (i.e., find a pair of commuting representatives or prove that none exist).
\item Find a full system of representatives.
\item Identify a full system of representatives so that the candidate maps are isomorphisms.
\eenum
For a right-angled Coxeter system, it happens that all of these steps are not just possible, but straightforward.

\begin{restatable}{Proposition}{invedges}
\label{prop:edge}
If $W_{\Gamma}$ is a right-angled Coxeter group, then two conjugacy classes of involutions $[x]$ and $[y]$ are connected by an edge in $\Inv[W_{\Gamma}]$ if and only if there exists another class $[z]$ such that $\overline{z} = \overline{xy}$ in the abelianization.
\end{restatable}
\pf
See Section \ref{sec:involutiongraph}.
\epf

If we are given a group $G$, supposing we can identify the conjugacy classes of involutions (i.e., the vertices of $\Inv[G]$), we can identify hypothetical edges and non-edges by looking for such $\overline{z}$ in $\Ab{G}$.  If $G$ is a right-angled Coxeter group, then this will produce the correct involution graph, and the remainder of the procedure will (hopefully, if we pick a good full system of representatives) identify a right-angled Coxeter presentation.  On the other hand, if this not-quite involution graph of $G$ is not a clique graph, we can be certain that $G$ is not a right-angled Coxeter group.  At no point do we directly need to check that we have calculated the true involution graph of $G$.  We summarize this discussion with the following amended collapsing procedure.  For remaining details (especially, a detailed description of how to do these calculations in the abelianization), refer to Section \ref{sec:involutiongraph}.

\begin{restatable}{Theorem}{algebraiccollapsing}
\label{thm:algebraiccollapsing}
Suppose $G$ is a group whose only torsion elements all have order 2 and so that $\Ab{G} \cong (\Z/2\Z)^n$ for some $n$.  If the following procedure returns True, then $G$ is a right-angled Coxeter group (and the procedure indicates a right-angled Coxeter presentation).  If the procedure returns False, then $G$ is not a right-angled Coxeter group.
\enum[labelindent=\parindent,label=\arabic*:]
\item Determine all conjugacy classes of involutions in $G$, and let these be the vertices of a graph $\Gamma'$.  If there are not finitely many, return \textbf{False}.
\item Apply Proposition \ref{prop:edge} to construct the edges of $\Gamma'$.
\item If $\Gamma'$ is not a clique graph, return \textbf{False}.
\item Find a full system of representatives for the vertices of $\Gamma'$.  If no such system exists, return \textbf{False}.
\item Collapse as in Theorem \ref{thm:collapsing}, using Proposition \ref{prop:choose} to ensure that nontrivial product relations are avoided.  Write $C(\Gamma')$ for the resulting graph.
\item Let $\Gamma$ be a graph isomorphic to $C(\Gamma')$ with generic vertex labels $a_1, \dotsc, a_n$.  Let $\varphi \colon W_{\Gamma} \to G$ be the map which sends the generators of $W_{\Gamma}$ to the word given by the corresponding labels of vertices in $C(\Gamma')$.  If $\varphi$ is an isomorphism, return \textbf{True}.
\item Otherwise, return \textbf{Unknown}.
\eenum
\end{restatable}

\section{Applications and Results}
\label{sec:applications}

In this section, we apply the recognition procedure from Section \ref{sec:summary} to seek out right-angled Coxeter presentations for certain families of groups.  We focus in particular on 
\enum
\item semi-direct products of a given right-angled Coxeter group $W_{\Gamma}$ by certain subgroups of $\Out^0(W_{\Gamma})$, and
\item subgroups of a given $W_{\Gamma}$ generated by chosen subsets of involutions.
\eenum
In particular, we note that the families of groups that we consider are already generated by involutions, have no torsion of order other than 2, and are usually given by presentations which are nearly right-angled Coxeter.  

If $D$ is a union of connected components of $\Gamma \setminus \St(a_i)$ for some $i$, then the automorphism of $W_{\Gamma}$ determined by
\[
\chi_{i, D}(a_j) = \begin{cases}
a_ia_ja_i & a_j\in D \\
a_j & \textrm{otherwise}
\end{cases}
\]
is called the \emph{partial conjugation with acting letter $a_i$ and domain $D$}. (Note that this terminology is not entirely consistent in the literature.  Other papers have reserved partial conjugation for the case in which $D$ is a single connected component \cite{GPR, CRSV}, while Laurence used the term \emph{locally inner automorphism} \cite{Laurence} before the term partial conjugation became common.  We have preferred here to allow for multiple connected components in the domain of a partial conjugation, and we would propose the term \emph{elementary partial conjugation} for the case in which $D$ consists of a single connected component.)  The partial conjugations generate $\Out^0(W_{\Gamma})$.

In Section \ref{sec:positiveresults}, we present families of groups which our procedure shows to be right-angled Coxeter.  One example is worked out in full detail to demonstrate the procedure.  For the remaining results, we simply state the resulting right-angled Coxeter group and the isomorphism determined by our procedure.  The reader is left to verify the details.  Most of these results are about split extensions of a given $W_{\Gamma}$ by a finite subgroup of $\Out^0(W_{\Gamma})$ generated by (pairwise commuting) partial conjugations.

In Section \ref{sec:negativeresults}, we present families of groups which our procedure shows cannot be right-angled Coxeter.  Again, one example is worked out in full detail.  We note one example which is of particular interest: we find a group $W_{\Gamma}$ with two elements $x, y\in \Out^0(W_{\Gamma})$ such that $G = W_{\Gamma} \rtimes \langle x, y\rangle$ is a right-angled Coxeter group, but $H = W_{\Gamma} \rtimes \langle xy \rangle$ is not.  In particular, we can realize $G$ as the iterated semi-direct product
\[
G = (W_{\Gamma} \rtimes \langle xy\rangle) \rtimes \langle x \rangle,
\]
where each extension has degree 2.  So this gives, to our knowledge, the first example in which the existence of a right-angled Coxeter presentation is lost and then recovered by semi-direct product extensions.

Finally, in Section \ref{sec:decompositionresults}, we note that many of our examples of right-angled Coxeter families arise as semi-direct products.  By analyzing the properties of the defining graphs of the groups arising from these semi-direct products, we can identify semi-direct product decompositions in many cases.  Such decompositions are generally not unique, and we cannot at the moment provide an exhaustive list of graph features of $\Gamma$ which indicate a semi-direct product decomposition of $W_{\Gamma}$.

\subsection{Groups Which are Right-angled Coxeter}
\label{sec:positiveresults}

\ex
\label{ex:racg}
We begin with an explicit example in which we demonstrate the recognition procedure in detail.  Consider the following defining graph:

\begin{figure}[h]
\begin{tikzpicture}[scale=1.5]

\node (a1) [circle,fill,inner sep=2pt,label=90:$a_1$] at (90:1) {};
\node (a2) [circle,fill,inner sep=2pt,label=210:$a_2$] at (210:1) {};
\node (a3) [circle,fill,inner sep=2pt,label=330:$a_3$] at (330:1) {};
\node (a4) [circle,fill,inner sep=2pt,label=270:$a_4$] at (0,0) {};

\draw (a1) to (a4);
\draw (a2) to (a4);
\draw (a3) to (a4);

\node at (-1.5,0.5) {$\Gamma$};

\end{tikzpicture}
\caption{The defining graph $\Gamma$.}
\label{fig:Gamma}
\end{figure}
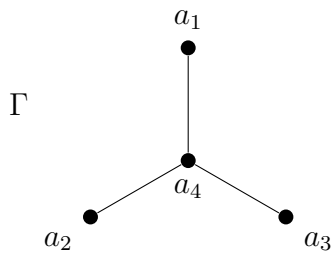

Write $x = \chi_{1, \{2\}}$ for the partial conjugation with acting letter $a_1$ and domain $\{a_2\}$. We consider the group $G = W_{\Gamma} \rtimes \langle x \rangle$, which has the following presentation:
\begin{align*}
G = \langle a_1, a_2, a_3, a_4, x &\mid a_i^2 = x^2 = 1, [a_1, a_4] = [a_2, a_4] = [a_3, a_4] = 1, \\
&\quad [a_1, x] = [a_3, x] = [a_4, x] = 1, xa_2x = a_1a_2a_1 \rangle
\end{align*}
This is not quite a right-angled Coxeter presentation, so we apply our procedure to see if we can find one.  

First, we compute $\Ab{G}$ (removing any relations that become trivial and understanding that group presentations with additive notation are assumed to be abelian): 
\begin{align*}
\Ab{G} &= \langle \overline{a_1}, \overline{a_2}, \overline{a_3}, \overline{a_4}, \overline{x} \mid 2 \overline{a_i} = 2 \overline{x} = 0\rangle \\
&\cong \langle \overline{a_1} \rangle \times \langle \overline{a_2} \rangle \times \langle \overline{a_3} \rangle \times \langle \overline{a_4} \rangle \times \langle \overline{x} \rangle \\
&\cong (\mathbb{Z}/2\mathbb{Z})^5
\end{align*}

The relation matrix
\[
\begin{pmatrix}
2 & 0 & 0 & 0 & 0 \\
0 & 2 & 0 & 0 & 0 \\
0 & 0 & 2 & 0 & 0 \\
0 & 0 & 0 & 2 & 0 \\
0 & 0 & 0 & 0 & 2
\end{pmatrix}
\]
is already in Smith normal form, and so our canonical abelianization map is just $G \to \Ab{G}: g \mapsto \overline{g}$.

We now want to list all conjugacy classes of involutions in $G$.  The classes of involutions in $W_{\Gamma}$ are evident by inspection of $\Gamma$: $a_i$ for each $i$, and $a_ja_4$ for each $1\leq j \leq 3$.  The new generator $x$ is also an involution, and the products of $x$ with the other generators that commute with it give new involutions: $xa_1, xa_3, xa_4$.  There are two remaining conjugacy classes of involutions, namely $xa_1a_2$ and $xa_1a_2a_4$. 

These are all of the conjugacy classes of involutions in $G$. We could try to prove this directly, but it will also end up following from the fact that our procedure in this case does in fact construct an explicit isomorphism with a right-angled Coxeter group. Thus, we can omit the details. 

We claim that the following is the involution graph $\Inv[G]$.  The given system of representatives is a full system, and the commuting relations are straightforward to check. (If they weren't as straightforward, we could easily construct the edge relations given by Proposition \ref{prop:edge}.)

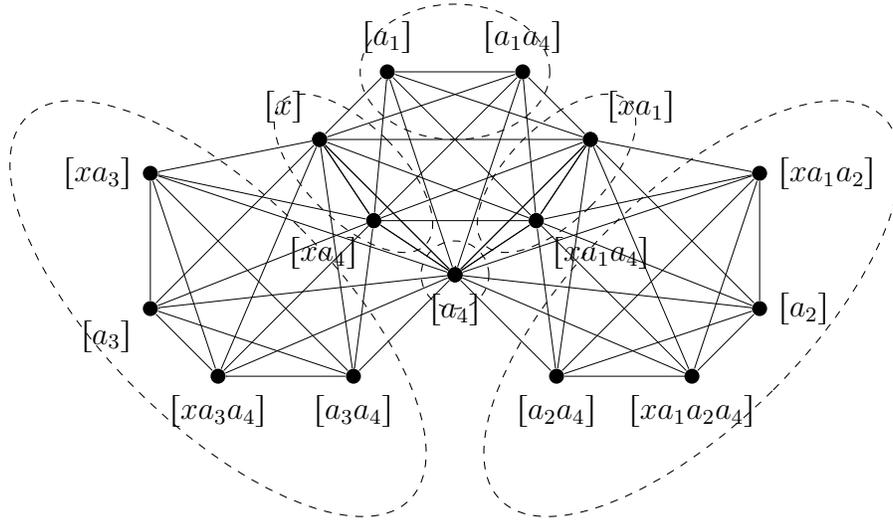
\begin{figure}[h]
\begin{tikzpicture}[scale=0.9]

\node (u1) [circle,fill,inner sep=2pt,label=180:{$[xa_3]$}] at (-4.5,1.5) {};
\node (u2) [circle,fill,inner sep=2pt,label=185:{$[a_3]$}] at (-4.5,-0.5) {};
\node (u3) [circle,fill,inner sep=2pt,label=-90:{$[xa_3a_4]$}] at (-3.5,-1.5) {};
\node (u4) [circle,fill,inner sep=2pt,label=-90:{$[a_3a_4]$}] at (-1.5,-1.5) {};

\node (u5) [circle,fill,inner sep=2pt,label=-90:{$[a_4]$}] at (0,0) {};
\node (u6) [circle,fill,inner sep=2pt,label=225:{$[xa_4]$}] at (-1.2,0.8) {};
\node (u7) [circle,fill,inner sep=2pt,label=135:{$[x]$}] at (-2,2) {};

\node (v1) [circle,fill,inner sep=2pt] at (0,0) {};
\node (v2) [circle,fill,inner sep=2pt] at (-1.2,0.8) {};
\node (v3) [circle,fill,inner sep=2pt] at (-2,2) {};

\node (v4) [circle,fill,inner sep=2pt,label=90:{$[a_1]$}] at (-1,3) {};
\node (v5) [circle,fill,inner sep=2pt,label=90:{$[a_1a_4]$}] at (1,3) {};

\node (v6) [circle,fill,inner sep=2pt,label=45:{$[xa_1]$}] at (2,2) {};
\node (v7) [circle,fill,inner sep=2pt,label=-45:{$[xa_1a_4]$}] at (1.2,0.8) {};

\node (w1) [circle,fill,inner sep=2pt] at (0,0) {};
\node (w2) [circle,fill,inner sep=2pt] at (1.2,0.8) {};
\node (w3) [circle,fill,inner sep=2pt] at (2,2) {};

\node (w4) [circle,fill,inner sep=2pt,label=0:{$[xa_1a_2]$}] at (4.5,1.5) {};
\node (w5) [circle,fill,inner sep=2pt,label=0:{$[a_2]$}] at (4.5,-0.5) {};
\node (w6) [circle,fill,inner sep=2pt,label=-90:{$[xa_1a_2a_4]$}] at (3.5,-1.5) {};
\node (w7) [circle,fill,inner sep=2pt,label=-90:{$[a_2a_4]$}] at (1.5,-1.5) {};

\foreach \x in {1,...,6}
  \foreach \y in {\x,...,7}
    {
      \draw [very thin] (u\x) to (u\y);
      \draw [very thin] (v\x) to (v\y);
      \draw [very thin] (w\x) to (w\y);
    }

\draw [dashed,rotate around={45:(-3.5,-0.5)}] (-3.5,-0.5) ellipse (1.7cm and 4cm);
\draw [dashed,rotate around={-45:(3.5,-0.5)}] (3.5,-0.5) ellipse (1.7cm and 4cm);
\draw [dashed,rotate around={45:(-1.5,1.5)}] (-1.5,1.5) ellipse (0.7cm and 1.5cm);
\draw [dashed,rotate around={-45:(1.5,1.5)}] (1.5,1.5) ellipse (0.7cm and 1.5cm);
\draw [dashed] (0,3) ellipse (1.4cm and 1cm);
\draw [dashed] (0,0) ellipse (0.5cm and 0.5cm);

\end{tikzpicture}
\caption{The involution graph $\Inv[G]$.}
\label{fig:involutiongraph}
\end{figure}

The brackets in the involution graph represent conjugacy classes.  Since we now have a full system of representatives, we may stop writing these brackets.  For the remainder of the calculation, brackets around a vertex label will denote its star equivalence class.  Before calculating the star poset structure, we observe that this graph clearly satisfies the Maximal Clique Condition and the Minimal Vertex Condition, and the Inclusion-Exclusion Condition is straightforward to verify.

The equivalence classes in the star poset are the following (identified by the dashed ellipses in the figure):
\begin{align*}
[a_1] &= \{a_1, a_1a_4\}                & [a_2]  &= \{a_2, a_2a_4, xa_1a_2, xa_1a_2a_4\} \\
[a_3] &= \{a_3, a_3a_4, xa_3, xa_3a_4\} & [a_4]  &= \{a_4\} \\
[x]   &= \{x, xa_4\}                    & [xa_1] &= \{xa_1, xa_1a_4\}
\end{align*}

The Hasse diagram for this poset is as follows:
\begin{figure}[h!]
\begin{tikzpicture}[scale=1.0]

\node (a4)  [circle,fill,inner sep=2pt,label=90:{$[a_4]$}]  at (2, 2) {};
\node (x)   [circle,fill,inner sep=2pt,label=135:{$[x]$}]   at (1, 1) {};
\node (xa1) [circle,fill,inner sep=2pt,label=45:{$[xa_1]$}]  at (3, 1) {};
\node (a3)  [circle,fill,inner sep=2pt,label=-90:{$[a_3]$}] at (0, 0) {};
\node (a1)  [circle,fill,inner sep=2pt,label=-90:{$[a_1]$}] at (2, 0) {};
\node (a2)  [circle,fill,inner sep=2pt,label=-90:{$[a_2]$}] at (4, 0) {};

\draw (a3) to (x) to (a4) to (xa1) to (a2);
\draw (x) to (a1) to (xa1);

\end{tikzpicture}
\caption{The Hasse diagram for the poset $\mathcal{P}(\Inv[G])$.}
\label{fig:Hasseexample}
\end{figure}
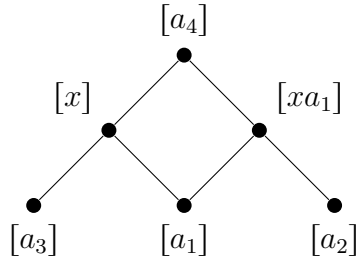

The element $[a_4]$ is maximal in the poset structure and contains a single element.  We add $a_4$ to $V$.  Next, we consider $[x]$ (or $[xa_1]$, the order in which we consider these classes is irrelevant).  The clique above $[x]$ has size 3, so 2 of its vertices must be added to $V$.  We have already added 1, so we must pick one more from $[x]$.  Examining the abelienization, $\langle \overline{a_4}, \overline{x} \rangle \cong (\mathbb{Z}/2\mathbb{Z})^2$ and either of $\overline{x}$ or $\overline{xa_4}$ will extend $\overline{a_4}$ into a basis. So we choose to add $x$ to $V$.  Similarly, we consider $[xa_1]$ and add $xa_1$ to $V$.

The remaining three classes are all minimal.  Suppose we take $[a_2]$ next.  The clique above $[a_2]$ has size 7, so we must choose 3 elements from it.  We have already chosen 2, so we need to choose 1 more. Checking the abelianization again, we see that any choice of the 4 elements in $[a_2]$ will extend to a basis, and so we add $a_2$ to $V$.  Similarly, from $[a_3]$, we add $a_3$ to $V$.

\begin{figure}[h]
\begin{tikzpicture}[scale=1.5]

\node (b1) [circle,fill,inner sep=2pt,label=-90:$a_3$] at (-1,0) {};
\node (b2) [circle,fill,inner sep=2pt,label=90:$x$] at (-0.5,0.87) {};
\node (b3) [circle,fill,inner sep=2pt,label=-90:$a_4$] at (0,0) {};
\node (b4) [circle,fill,inner sep=2pt,label=90:$xa_1$] at (0.5,0.87) {};
\node (b5) [circle,fill,inner sep=2pt,label=-90:$a_2$] at (1,0) {};

\draw (b1) to (b2) to (b3) to (b4) to (b5) to (b3) to (b1);
\draw (b2) to (b4);

\begin{scope}[xshift=3.5cm]
\node (b1) [circle,fill,inner sep=2pt,label=-90:$b_1$] at (-1,0) {};
\node (b2) [circle,fill,inner sep=2pt,label=90:$b_2$] at (-0.5,0.87) {};
\node (b3) [circle,fill,inner sep=2pt,label=-90:$b_3$] at (0,0) {};
\node (b4) [circle,fill,inner sep=2pt,label=90:$b_4$] at (0.5,0.87) {};
\node (b5) [circle,fill,inner sep=2pt,label=-90:$b_5$] at (1,0) {};

\draw (b1) to (b2) to (b3) to (b4) to (b5) to (b3) to (b1);
\draw (b2) to (b4);
\end{scope}

\end{tikzpicture}
\caption{On the left is the collapsed graph $\Lambda$.  On the right is an isomorphic graph with generic labels.}
\label{fig:Lambda}
\end{figure}
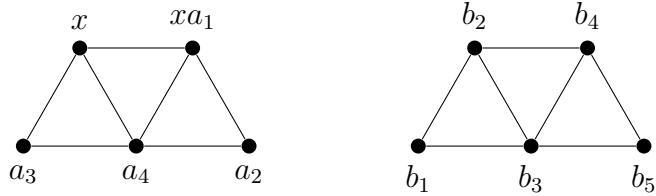

Finally, we consider $[a_1]$.  The clique above $[a_1]$ has size 7, and we have already chosen 3 of these vertices, so we choose no more.  This leaves us with $V = \{a_2, a_3, a_4, x, xa_1\}$.  We take the induced subgraph $\Lambda$ of $\Inv[G]$ on these vertices (Figure \ref{fig:Lambda}).

We now have a candidate map $\varphi \colon W_{\Lambda} \to G$.  It is straightforward to check that the map $\psi$ below is the inverse, and that $\varphi$ and $\psi$ are isomorphisms:
\begin{align*}
\varphi\colon b_1 &\mapsto a_3 & \psi\colon a_1 &\mapsto b_2 b_4 \\
b_2 &\mapsto x & a_2 &\mapsto b_5 \\
b_3 &\mapsto a_4 & a_3 &\mapsto b_1 \\
b_4 &\mapsto xa_1 & a_4 &\mapsto b_3 \\
b_5 &\mapsto a_2 & x &\mapsto b_2 
\end{align*}
Thus, $G$ is a right-angled Coxeter group, completing the example.
\eex

In this example, we were extending a right-angled Coxeter group by a single partial conjugation.  It turns out that this will always yield a right-angled Coxeter group, and in fact we can say much more.

\lem
\label{lem:pcextensions}
Suppose $W_\Gamma$ is a right-angled Coxeter group.  If $\alpha_1, \dots, \alpha_k$ are partial conjugations of $W$ with the same acting letter and pairwise disjoint domains, then $G = W \rtimes \langle \alpha_1, \dots, \alpha_k \rangle$ is a right-angled Coxeter group.
\elem
\pf
Without loss of generality, we may assume each $\alpha_j$ has acting letter $a_1$.  Let $D_i$ denote the domain of $\alpha_i$ for each $1 \leq 1 \leq k$.  Now $G$ is generated by the elements $\{a_1, \dotsc, a_n, \alpha_1, \dotsc, \alpha_k\}$ and presented by the following relations:
\enum[labelindent=\parindent,label=(R\arabic*)]
\item $a_i^2=1$ for $1\leq i \leq n$,
\item $[a_i, a_j]=1$ for $\{a_i, a_j\} \in E(\Gamma)$,
\item $\alpha_i^2=1$ for $1 \leq i \leq k$,
\item $[\alpha_i, \alpha_j]=1$ for $1 \leq i < j \leq k$ 
\item $[\alpha_i, a_j]=1$ for $a_j\notin D_j$, and
\item $\alpha_i a_j \alpha_i = a_1 a_j a_1$ for $a_j\in D_i$.
\eenum

Let $H$ be the group generated by $\{b_1, \dotsc, b_n, \beta_1, \dots, \beta_k\}$ and presented by the relations:
\enum[labelindent=\parindent,label=(S\arabic*)]
\item $b_i^2=1$ for $1\leq i \leq n$,
\item $[b_i, b_j]=1$ for $\{a_i, a_j\} \in E(\Gamma)$,
\item $\beta_i^2=1$ for $1 \leq i \leq k$,
\item $[\beta_i, \beta_j]=1$ for $1 \leq i < j \leq k$ 
\item $[\beta_i, b_j]=1$ for $a_j\notin D_i$, and
\item $[b_1, b_i]=1$ for $2 \leq i \leq n$ and $a_i \in D_1 \cup \dots \cup D_k$.
\eenum
We note that the given presentation for $H$ is a right-angled Coxeter presentation.  We define maps
\begin{align*}
\widehat{\varphi} \colon \{a_1, \dotsc, a_n, \alpha_1, \dotsc, \alpha_k\} &\to \{b_1, \dotsc, b_n, \beta_1, \dotsc, \beta_k\} \\
a_1 &\mapsto b_1 \beta_1 \dotsc \beta_k \\
\alpha_i &\mapsto \beta_i \ (1\leq i \leq k) \\
a_i &\mapsto b_i\ (2\leq i \leq n) \\ 
\widehat{\psi} \colon \{b_1, \dotsc, b_n, \beta_1, \dotsc, \beta_k\} &\to \{a_1, \dotsc, a_n, \alpha_1, \dotsc, \alpha_k\} \\
b_1 &\mapsto a_1 \alpha_1 \dots \alpha_k \\
\beta_i &\mapsto \alpha_i\ (1\leq i \leq k)\\
b_i &\mapsto a_i\ (2\leq i \leq n)
\end{align*}
It is straightforward to check that $\widehat{\varphi}$ and $\widehat{\psi}$ preserve the relations (R1)-(R6) and (S1)-(S6), respectively, so they induce homomorphisms $\varphi \colon G \to H$ and $\psi \colon H \to G$.  (Note that the preservation of the relation (S6) uses the assumption that the domains $D_i$ are pairwise disjoint.)  Finally, it is straightforward to see that $\varphi$ and $\psi$ are inverses to each other, hence $G$ and $H$ are isomorphic.  That is, $G$ is a right-angled Coxeter group.
\epf

Suppose $H \leq \Out^0(W_{\Gamma})$ is generated by partial conjugations $\chi_1, \dotsc, \chi_k$.  Having shown that the semi-direct product extension of $W_{\Gamma}$ by any single partial conjugation is again right-angled Coxeter, we might hope to show that $W_{\Gamma} \rtimes H$ is right-angled Coxeter by observing that this is isomorphic to taking the iterated semi-direct products, each by a single $\chi_i$:
\[
W_{\Gamma} \rtimes H = \left(\dotsb \left(\left(W_{\Gamma} \rtimes \langle \chi_1\rangle\right) \rtimes \langle \chi_2 \rangle\right) \rtimes \dotsb \rtimes \langle \chi_k\rangle \right).
\]
However, there is a subtlety that ruins this argument, namely, that $\chi_2$ will extend to some automorphism of $W_{\Gamma} \rtimes \langle \chi_1 \rangle$, but not necessarily to a partial conjugation.  We cannot extend inductively, since we cannot ensure that we are always extending by single partial conjugations.  The following lemma and theorem identifies certain cases in which this inductive argument works.

\lem
\label{lem:gammaisPC}
Suppose $W, \Gamma, a_1, \alpha_1, \dots, \alpha_k, H$ and $G$ are as in the lemma and proof above.  Let $\gamma$ be a partial conjugation of $W$ with acting letter $a_2 \neq a_1$ and such that $\gamma$ commutes with each of the automorphisms $\alpha_1, \dots, \alpha_k$.  Then $\gamma$ acts on $G$ as a partial conjugation.
\elem
\pf
Without loss of generality we may assume $\gamma$ has acting letter $a_2$ and domain $D$.  Recall that $a_2 = b_2$.  To show that $\gamma$ acts on $G$ as a partial conjugation we shall consider the result of conjugation by $\gamma$ on each of the generators $b_1, \dots, b_n, \beta_1, \dots, \beta_k$. Firstly we note: $\gamma \beta_i \gamma = \beta_i$ for $1 \leq i \leq k$; $\gamma b_i \gamma = b_i$ for $1 \leq i \leq n$ and $a_i \notin D$; $\gamma b_i \gamma = b_2 b_i b_2$ for $2 \leq i \leq n$ and $a_i \in D$.  If $a_1 \notin D$, then $\gamma b_1 \gamma = \gamma a_1 \gamma = b_1$.  Suppose  $a_1 \in D$.  Since $\gamma$ commutes pairwise with $\alpha_1, \dots, \alpha_k$, we have that $a_2 \notin D_1 \cup \dots \cup D_k$.   We compute
\begin{align*}
\gamma b_1 \gamma &= \gamma a_1 \alpha_1 \dotsc \alpha_k \gamma \\
&= \gamma a_1 \gamma \alpha_1 \dotsc \alpha_k \\
&= a_2 a_1 a_2\alpha_1 \dots \alpha_k \\
&= a_2 a_1 \alpha_1 \dots \alpha_k a_2 \\
&= b_2 b_1 b_2.
\end{align*}
Since $\gamma$ is an automorphism of $G$, and $\gamma$ takes each generator to either itself, or the conjugate of itself by $b_2$, $\gamma$ is a partial conjugation of $G$.

Write $\varphi \colon \{a_1, \dotsc, a_n\} \to \{b_1, \dotsc, b_n\}$ for the map $\varphi(a_i) = b_i$.  From the calculations above, the domain of $\gamma$ acting on $G$ is $\varphi(D)$.
\epf

\ExtensionByPCs
\pf
The proof is by induction, applying the lemmas above at each step.  Let $\alpha_1, \dotsc, \alpha_{k_1}$ be those $\chi_i$ with acting letter 1.  By assumption, they have pairwise disjoint domains.  By Lemma \ref{lem:pcextensions}, $W_{\Gamma} \rtimes \langle \alpha_1, \dotsc, \alpha_{k_1}\rangle$ is a RACG.

Moreover, by Lemma \ref{lem:gammaisPC}, the remaining $\chi_i$'s still act like partial conjugations, and their domains do not intersect, since they didn't before the extension.  Now take $\beta_1, \dotsc, \beta_{k_2}$ among the remaining $\chi_i$ to be those which have acting letter 2, and extend by $\langle \beta_1, \dotsc, \beta_{k_2}\rangle$.

Continuing inductively, we extend at the $i$th step by all remaining partial conjugations with acting letter $i$.  The result follows.
\epf

In \cite{GPR}, the authors investigate the automorphism groups of graph products of cyclic groups. In the case that $W$ is a right-angled Coxeter group, the authors recover a result from \cite{Tits} which shows $\Aut(W) = \Aut^0(W)\rtimes \Aut^1(W)$ with $\Aut^1(W)$ finite. Thus $\Aut^0(W)$ (sometimes denoted $\Aut^{PC}(W)$) which is the subgroup of $\Aut(W)$ generated by all partial conjugations of $W$, is a finite index subgroup of $\Aut(W)$. Furthermore, they show that $\Aut^0(W)$ splits as $\Inn(W) \rtimes \Out^0(W)$. Finally, they give the following condition on $\Gamma$, called \emph{no SILs}, which  characterizes exactly when $\Out^0(W)$ is finite and is thus isomorphic to $\mathbb{Z}_2^n$.

\defi
A graph $\Gamma$ has a \emph{separating intersection of links (SIL)} if, for some vertices $v$ and $w$ with $d(v, w) \geq 2$, there is a component of $\Gamma \setminus (\Lk(v) \cap \Lk(w))$ which contains neither $v$ nor $w$.  Otherwise, $\Gamma$ is said to have \emph{no SILs}.
\edefi

$\Inn(W_{\Gamma})$ is known to be a right-angled Coxeter group.  In the case that $\Gamma$ has no SILs, $\Aut^0(W_{\Gamma})$ is a finite extension of $\Inn(W_{\Gamma})$. In \cite{CRSV}, it is shown that $\Aut^0(W_{\Gamma})$ is again a right-angled Coxeter group in that case.  We arrive at this same result as a direct application of the previous corollary.

\cor
If $\Gamma$ contains no SILs, then $\Aut^0(W)$ is a right-angled Coxeter group and thus $\Aut(W)$ contains a right-angled Coxeter group as a subgroup of finite index. 
\ecor
\pf
Without loss of generality we may assume $W$ has trivial center.  Suppose $\Gamma$ contains no SILs.  Then
\[
\Aut^0(W) = \Inn(W) \rtimes \Out^0(W) \cong W \rtimes \Out^0(W),
\]
and $\Out^0(W)$ is generated by pairwise commuting partial conjugations which satisfy the condition in the corollary above.  
\epf

In general, one should not expect $\Aut(W)$ to be right-angled Coxeter.  The elements of $\Aut^1(W)$ include graph symmetries, which could then introduce torsion elements of order other than 2.  One should not generally expect that $\Aut^0(W)$ is a right-angled Coxeter group, but one might see the no SILs result as suggesting that we restrict our attention to extensions of right-angled Coxeter groups by finite subgroups of $\Out^0(W)$ (although Example \ref{ex:extofext} in the following section demonstrates that even this restriction is not sufficient).

\subsection{Groups Which are not Right-angled Coxeter}
\label{sec:negativeresults}

\ex
\label{ex:nonracg}
As in the previous section, we begin with an explicitly worked out example. Let $G$ denote the group presented as follows:
\begin{align*}
G = \langle a, b, c, x, y &\mid a^2, b^2, c^2, x^2, y^2, \\
&\quad xax = a, xbx = b, xcx = aca, \\
&\quad yay=a, yby=b, ycy=bcb\rangle.
\end{align*}
Let $W = \langle a, b, c \rangle$ and $H = \langle x, y \rangle$.  Then $W = \mathbb{Z} / 2\mathbb{Z} \ast \mathbb{Z} / 2\mathbb{Z} \ast \mathbb{Z} / 2\mathbb{Z}$, $H \cong \mathbb{Z} / 2\mathbb{Z} \ast \mathbb{Z} / 2\mathbb{Z}$ and $G \cong W \rtimes H$, where $x$ and $y$ act as a pair of non-commuting partial conjugations.

To construct $\Inv[G]$, we must understand the involutions in $G$.  Since $G = W \rtimes H$, each $g \in G$ may be written uniquely in the form $g = wh$, where $w \in W$ and $h \in H$.  Further, $g^2 = whwh = whw(h^{-1}h)h = ww^{h^{-1}} h^2$. Since every element in $G$ can be uniquely written as a product of an element of $G$ and an element of $H$, if $g$ is an involution, then $h$ is an involution and $w^{h^{-1}} = w^h = w^{-1}$. $H$ is a right-angled Coxeter group (in fact, $D_{\infty}$), and so every non-trivial involution in $H$ is conjugate to either $x$ or $y$; it follows that, up to conjugation, we may suppose $g$ has one of the following forms:
\enum
\item $w$ such that $w^2 = 1$,
\item $wx$ such that $w^x = w^{-1}$, or
\item $wy$ such that $w^y = w^{-1}$.
\eenum
Every element of the first type is conjugate to either $a$, $b$, or $c$.  Now we'll try to list elements of the second type (elements of the third type will be analogous).

Suppose $g = wx$ with $w^x = w^{-1}$.  We further suppose that, within the collection of words of this form in the conjugacy class of $g$, we choose the shortest possible $w$.  The element $w$ can be written uniquely in the form $u_0 b u_1 b \dotsm u_{m-1} b u_m$ where $m \geq 0$, each $u_i$ is a geodesic word in $\{a, c\}^\ast$, and only $u_0$ and $u_m$ may be trivial.  Then $w^x = u_0^x b u_1^x b \dotsm u_{m-1}^x b u_m^x = w^{-1}$ implies that $u_0^x = u_m^{-1}, u_1^x = u_{m-1}^{-1}$, and so on.  We now consider a few subcases.

If $m > 0$ and $u_0$ is not trivial, then
\begin{align*}
u_0^{-1} (wx) u_0 &= u_0^{-1} (u_0 b u_1 b \dotsm u_{m-1} b u_m x) u_0 \\
&= b u_1 b \dotsm u_{m-1} b u_m u_0^x x \\
&= b u_1 b \dotsm u_{m-1} b x.
\end{align*}
This contradicts the minimality of the length of $w$, so either $m = 0$ or $u_0$ is trivial.  If $u_0$ is trivial and $m > 1$, then $w$ begins and ends with $b$, in which case $|b(wx)b| < |wx|$.  Again, this contradicts minimality, hence either $m=0$ or $w = b$.

If $m = 0$, then $w = u_0 \in \langle a, c \rangle$ is geodesic and so is an alternating string of $a$ and $c$.  If $|w| > 1$ and $|w|$ is odd, then $w$ begins and ends with the same letter.  If $w$ begins and ends with $a$, then $|a wx a| = |a w a x | < |wx|$; if $w$ begins and ends with $c$ then $w^x$ begins and ends with $a$, hence $w^x \neq w^{-1}$.  In either case, we have a contradiction, so $|w| = 1$, in which case $w=a$ or $w=c$, or else $|w|$ is even.  If $w = (ac)^n$ and $n > 1$, then $|aca(wx)aca| < |wx|$; if $w = (ca)^n$ and $n > 1$, then $|cwxc| < |wx|$.  In both cases, we have a contradiction. Our only case left is $m=0, n=1$, which corresponds to $w = ac$ or $w=ca$. Therefore, our only non-trivial possibilities for $w$ are $w = b, a, c, ac, ca$.

Note that $a(cax)a = acx$, so these cases fall into the same conjugacy classes.  In summary, we have that each involution of the form $wx$ is conjugate to exactly one of the elements $x, ax, bx, acx$.  (We observe that the final option $cx$ is not, in fact, an involution.  In this case, $w = c$, and $w^x \neq w^{-1}$.) We also observe that none of these involutions are conjugate to each other since they all map to distinct elements in $\Ab{G}$. 

Similarly, each involution of the form $wy$ is conjugate to exactly one of the elements $y, ay, by, bcy$.  Therefore the following is the complete list of conjugacy classes in $G$, and hence serves as the list of vertex labels in $\Inv[G]$:
\[
[a], [b], [c], [x], [ax], [bx], [acx], [y], [ay], [by], [bcy].
\]

\begin{figure}[ht]
\begin{tikzpicture}[xscale=1.8,yscale=1.2]

\node (acx) [circle,fill,inner sep=2pt,label=-90:{$[acx]$}] at (-2,0) {};
\node (c)   [circle,fill,inner sep=2pt,label=-135:{$[c]$}] at (0,0) {};
\node (bcy) [circle,fill,inner sep=2pt,label=-90:{$[bcy]$}] at (2,0) {};
\node (ax)  [circle,fill,inner sep=2pt,label=135:{$[ax]$}] at (-2,1.2) {};
\node (a)   [circle,fill,inner sep=2pt,label=180:{$[a]$}] at (-2,3) {};
\node (x)   [circle,fill,inner sep=2pt,label=-90:{$[x]$}] at (-1.5,1.5) {};
\node (bx)  [circle,fill,inner sep=2pt,label=-90:{$[bx]$}] at (0,0.8) {};
\node (b)   [circle,fill,inner sep=2pt,label=-90:{$[b]$}] at (1.5,1.5) {};
\node (by)  [circle,fill,inner sep=2pt,label=45:{$[by]$}] at (2,1.2) {};
\node (y)   [circle,fill,inner sep=2pt,label=0:{$[y]$}] at (2,3) {};
\node (ay)  [circle,fill,inner sep=2pt,label=90:{$[ay]$}] at (0,3.6) {};

\draw (a) to (x);
\draw (a) to (ax);
\draw (a) to (y);
\draw (a) to (ay);
\draw (b) to (x);
\draw (b) to (bx);
\draw (b) to (y);
\draw (b) to (by);
\draw (c) to (ax);
\draw (c) to (acx);
\draw (c) to (by);
\draw (c) to (bcy);
\draw (x) to (ax);
\draw (x) to (bx);
\draw (ax) to (acx);
\draw (y) to (ay);
\draw (y) to (by);
\draw (by) to (bcy);

\node at (-1.5,0.4) {$\Gamma_1$};
\node at (1.5,0.4) {$\Gamma_2$};
\node at (-1.8,1.7) {$\Gamma_3$};
\node at (1.8,1.7) {$\Gamma_4$};
\node at (0,1.2) {$\Gamma_5$};
\node at (0,3.2) {$\Gamma_6$};

\end{tikzpicture}
\caption{An involution graph which cannot be a clique graph.  The labeled triangles $\Gamma_i$ are the maximal cliques.}
\label{fig:nonexample}
\end{figure}
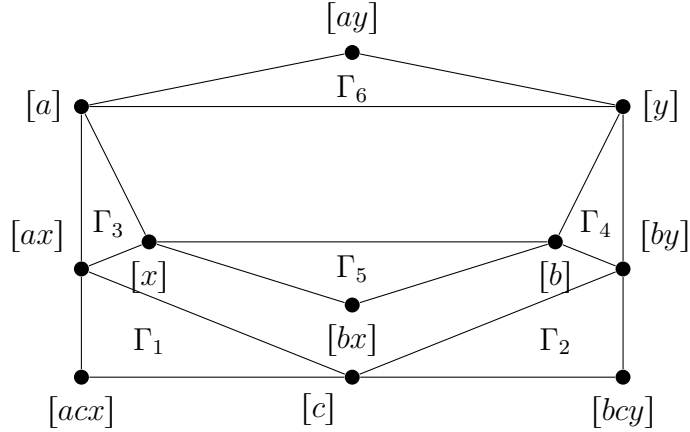

We now consider pairs of distinct conjugacy classes, to see whether or not they should be adjacent in $\Inv[G]$. By Proposition \ref{prop:edge}, we can just check the product relations among the images of the involutions in $\Ab{G}$. We omit the actual calculation and show the resulting involution graph in Figure \ref{fig:nonexample}.

Now $\Inv[G]$ is not a clique graph, since, for example, the Inclusion-Exclusion Condition fails.  (The reader can check this directly for the maximal cliques labeled $\Gamma_3$ and $\Gamma_4$ in the figure.)  This completes the example.
\eex

\begin{figure}[ht]
\begin{tikzpicture}[xscale=1.8,yscale=1.2]

\node (acx) [circle,fill,inner sep=2pt,label=-90:{$[acx]$}] at (-2,0) {};
\node (c)   [circle,fill,inner sep=2pt,label=-135:{$[c]$}] at (0,0) {};
\node (bcy) [circle,fill,inner sep=2pt,label=-90:{$[bcy]$}] at (2,0) {};
\node (ax)  [circle,fill,inner sep=2pt,label=135:{$[ax]$}] at (-2,1.2) {};
\node (a)   [circle,fill,inner sep=2pt,label=180:{$[a]$}] at (-2,3) {};
\node (x)   [circle,fill,inner sep=2pt,label=-90:{$[x]$}] at (-1.5,1.5) {};
\node (bx)  [circle,fill,inner sep=2pt,label=-90:{$[bx]$}] at (0,0.8) {};
\node (b)   [circle,fill,inner sep=2pt,label=-90:{$[b]$}] at (1.5,1.5) {};
\node (by)  [circle,fill,inner sep=2pt,label=45:{$[by]$}] at (2,1.2) {};
\node (y)   [circle,fill,inner sep=2pt,label=0:{$[y]$}] at (2,3) {};
\node (ay)  [circle,fill,inner sep=2pt,label=90:{$[ay]$}] at (0,3.6) {};

\node (z)  [circle,fill,inner sep=2pt,label=180:{$[z]$}] at (-1.2, 2) {};
\node (az)  [circle,fill,inner sep=2pt,label=0:{$[az]$}] at (-1.2,2.7) {};
\node (cz)  [circle,fill,inner sep=2pt,label=135:{$[cz]$}] at (1.2,2) {};
\node (cbz)  [circle,fill,inner sep=2pt,label=180:{$[cbz]$}] at (1.2,2.7) {};

\draw (a) to (x);
\draw (a) to (ax);
\draw (a) to (y);
\draw (a) to (ay);
\draw (b) to (x);
\draw (b) to (bx);
\draw (b) to (y);
\draw (b) to (by);
\draw (c) to (ax);
\draw (c) to (acx);
\draw (c) to (by);
\draw (c) to (bcy);
\draw (x) to (ax);
\draw (x) to (bx);
\draw (ax) to (acx);
\draw (y) to (ay);
\draw (y) to (by);
\draw (by) to (bcy);

\draw (a) to (az) to (z) to (a);
\draw (b) to (cz) to (cbz) to (b);
\draw (z) to (c) to (cz) to (z);


\end{tikzpicture}
\caption{The involution graph for $\Aut^0(W_3)$.}
\label{fig:Aut0W3}
\end{figure}
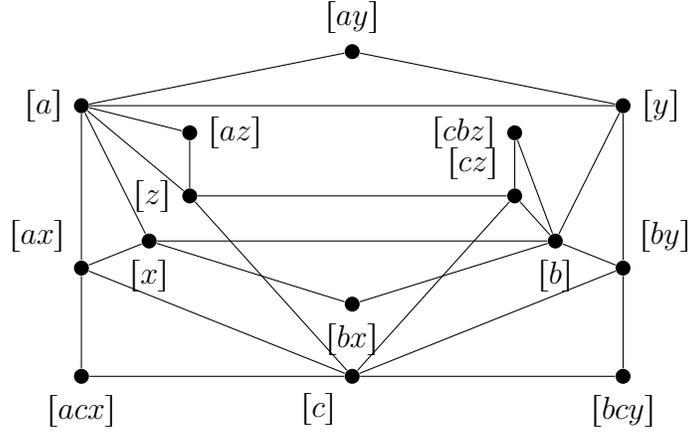

\ex
$\Aut^0(W_3)$ is not a right-angled Coxeter group.  The details are very similar to the previous example (we extend by one further partial conjugation), and are omitted here.  The involution graph is shown in Figure \ref{fig:Aut0W3}.

\eex

Here we must give the following warning.  The proof above relies on finding a portion of the involution graph which we know should not appear in any clique graph.  In the example, it is the ``triangle of triangles'' configuration (see Example \ref{ex:triangleoftriangles}).  This should not occur in the involution graph of a right-angled Coxeter group, essentially because it means that all three vertices of the central triangle must be generators (whereas, by construction of the involution graph in the case of right-angled Coxeter groups, we should expect two of the vertices are generators and the third is their product).

However, we must point out that, strictly speaking, there is no such thing as a ``poison pill'' subgraph---a subgraph which, by its presence, prevents the given graph from being a clique graph.  Indeed, if $\Gamma$ is any graph, then $\Gamma$ is an induced subgraph of $\Gamma_K$.  In this way, any finite graph may appear as an induced subgraph in some clique graph (even the ``triangle of triangles'').  In the example above, it is important that we know the central triangles $\Gamma_3$ and $\Gamma_4$ to be not just induced subgraphs, but also maximal cliques.

In all of the previous results, we have only considered split extensions by subgroups $H \leq \Out^0(W_{\Gamma})$ which were generated by partial conjugations.  In particular, if the partial conjugations commuted pairwise, then $H$ was finite and the extension $G = W_{\Gamma} \rtimes H$ was right-angled Coxeter.  On the other hand, in the example above, the partial conjugations did not commute, thus $H$ was infinite and $G$ was not right-angled Coxeter.  One might wonder whether the existence of a right-angled Coxeter presentation for the extension $G$ depends only on the finiteness of $H$.  The following example answers this question in the negative.

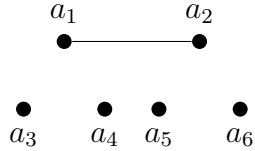
\begin{figure}[h]
\begin{tikzpicture}[scale=1.8]

\node (a1) [circle,fill,inner sep=2pt,label=90:$a_1$] at (0,0.5) {};
\node (a2) [circle,fill,inner sep=2pt,label=90:$a_2$] at (1,0.5) {};
\node (a3) [circle,fill,inner sep=2pt,label=-90:$a_3$] at (-0.3,0) {};
\node (a4) [circle,fill,inner sep=2pt,label=-90:$a_4$] at (0.3,0) {};
\node (a5) [circle,fill,inner sep=2pt,label=-90:$a_5$] at (0.7,0) {};
\node (a6) [circle,fill,inner sep=2pt,label=-90:$a_6$] at (1.3,0) {};

\draw (a1) to (a2);

\end{tikzpicture}
\caption{The defining graph $\Gamma$.}
\label{fig:extofext}
\end{figure}

\ex
\label{ex:extofext}
Let $\Gamma$ be the graph shown in Figure \ref{fig:extofext}.

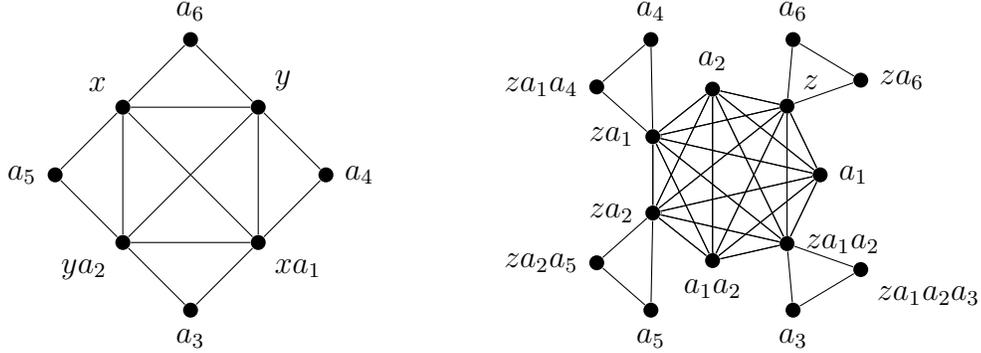
\begin{figure}[h]
\begin{tikzpicture}[scale=0.9]

\node (a1) [circle,fill,inner sep=2pt,label=-45:$xa_1$] at (1,-1) {};
\node (a2) [circle,fill,inner sep=2pt,label=225:$ya_2$] at (-1,-1) {};
\node (a3) [circle,fill,inner sep=2pt,label=-90:$a_3$] at (0,-2) {};
\node (a4) [circle,fill,inner sep=2pt,label=0:$a_4$] at (2,0) {};
\node (a5) [circle,fill,inner sep=2pt,label=180:$a_5$] at (-2,0) {};
\node (a6) [circle,fill,inner sep=2pt,label=90:$a_6$] at (0,2) {};
\node (x) [circle,fill,inner sep=2pt,label=135:$x$] at (-1,1) {};
\node (y) [circle,fill,inner sep=2pt,label=45:$y$] at (1,1) {};

\draw (a1) to (a2) to (x) to (y) to (a1) to (x);
\draw (a2) to (y) to (a4) to (a1) to (a3) to (a2) to (a5) to (x) to (a6) to (y);

\begin{scope}[xshift=8cm]
\node (v1) [circle,fill,inner sep=2pt,label={0*(360/7)}:$a_1$] at ({0*(360/7)}:1.3) {};
\node (v2) [circle,fill,inner sep=2pt,label={1*(360/7)}:$z$] at ({1*(360/7)}:1.3) {};
\node (v3) [circle,fill,inner sep=2pt,label=90:$a_2$] at ({2*(360/7)}:1.3) {};
\node (v4) [circle,fill,inner sep=2pt,label=180:$za_1$] at ({3*(360/7)}:1.3) {};
\node (v5) [circle,fill,inner sep=2pt,label=180:$za_2$] at ({4*(360/7)}:1.3) {};
\node (v6) [circle,fill,inner sep=2pt,label=-90:$a_1a_2$] at ({5*(360/7)}:1.3) {};
\node (v7) [circle,fill,inner sep=2pt,label=0:$za_1a_2$] at ({6*(360/7)}:1.3) {};

\node (a6) [circle,fill,inner sep=2pt,label=90:$a_6$] at (0.9,2) {};
\node (za6) [circle,fill,inner sep=2pt,label=0:$za_6$] at (1.9,1.4) {};

\node (a3) [circle,fill,inner sep=2pt,label=-90:$a_3$] at (0.9,-2) {};
\node (za1a2a3) [circle,fill,inner sep=2pt,label=-45:$za_1a_2a_3$] at (1.9,-1.4) {};

\node (a4) [circle,fill,inner sep=2pt,label=90:$a_4$] at (-1.2,2) {};
\node (za1a4) [circle,fill,inner sep=2pt,label=180:$za_1a_4$] at (-2,1.3) {};

\node (a5) [circle,fill,inner sep=2pt,label=-90:$a_5$] at (-1.2,-2) {};
\node (za2a5) [circle,fill,inner sep=2pt,label=180:$za_2a_5$] at (-2,-1.3) {};

\foreach \x in {1,2,...,7}
  \foreach \y in {1,2,...,7}
    {
      \draw (v\x) to (v\y);
    }

\draw (v2) to (a6) to (za6) to (v2);
\draw (v7) to (a3) to (za1a2a3) to (v7);
\draw (v4) to (a4) to (za1a4) to (v4);
\draw (v5) to (a5) to (za2a5) to (v5);

\end{scope}

\end{tikzpicture}
\caption{The defining graph of $G$ (left) and the involution graph of $H$ (right).}
\label{fig:graphofbigext}
\end{figure}

Let $x$ be the partial conjugation with acting letter $a_1$ and domain $\{a_3, a_4\}$, and let $y$ be the partial conjugation with acting letter $a_2$ and domain $\{a_3, a_5\}$.  Since $a_1$ and $a_2$ commute, so do $x$ and $y$.  Write $z = xy$ for the product, which is also an involution.  It follows from Theorem \ref{thm:mainextensiontheorem} that $G = W_{\Gamma} \rtimes \langle x, y \rangle$ is a right-angled Coxeter group.  Consider the subgroup $H = W_{\Gamma} \rtimes \langle z \rangle \leq G$.  The defining graph for $G$ and the involution graph for $H$ are shown in Figure \ref{fig:graphofbigext}.

The reader could verify $\Inv[H]$ in two ways---first, by directly calculating the involutions and checking their commuting relations; and second, using the defining graph of $G$ to calculate $\Inv[G]$, and then picking out the subset of vertices in $\Inv[G]$ which are labeled by elements in the subgroup $H$.  (Note that this latter method of constructing the involution graph of a subgroup will not work in general.  It works for the current example because $G$ is a right-angled Coxeter group and $H$ is normal.)

We can realize $G$ as the iterated semi-direct product
\[
G = \left(W_{\Gamma} \rtimes \langle z \rangle\right) \rtimes \langle x \rangle = H \rtimes \langle x \rangle.
\]
This gives an example of a right-angled Coxeter group $W_{\Gamma}$ with a degree 2 split extension $H$ which is not right-angled Coxeter.  Moreover, taking a further degree 2 extension $G$, we recover right-angled Coxeterness.
\eex

\subsection{Semi-direct Product Decompositions}
\label{sec:decompositionresults}

Here we present some results which are unrelated to the problem of recognizing right-angled Coxeter groups.  These results fall naturally out of the applications in Section \ref{sec:positiveresults}, and they generally address our ability to recognize semi-direct product decompositions of $W_{\Gamma}$ by identifying features of $\Gamma$.

To give the basic idea of how to generate these results, we give the following alternate description of Lemma \ref{lem:pcextensions}.  Suppose $a_1, \dotsc, a_n$ are the vertices of $\Gamma$ and $\alpha_1, \dotsc, \alpha_k$ are partial conjugations as in the lemma.  We will suppose that $a_1$ is the acting letter and $D_i$ is the domain of $\alpha_i$.  The lemma says that the group $G = W_{\Gamma} \rtimes \langle \alpha_1,\dotsc, \alpha_k \rangle$ is a right-angled Coxeter group, and the proof of the lemma gives the right-angled Coxeter generating set.  We can directly construct the defining graph $\Lambda$ for $G$ from $\Gamma$ as follows:
\enum
\item Add $k$ new vertices labeled $\alpha_1, \dotsc, \alpha_k$, all connected to one another and to $a_1$.
\item Connect each $\alpha_i$ to every $a_j$ where $a_j \notin D_i$.
\item Relabel $a_1$ as $a_1\alpha_1\alpha_2\dotsm \alpha_k$, and connect this to each vertex in $D_1 \cup D_2 \cup \dotsb \cup D_k$.
\eenum
The vertices $a_1, \alpha_1, \dotsc, \alpha_k$ form a clique of size $k+1$, and the union of the stars of these vertices cover all of $\Lambda$.  The restriction in Lemma \ref{lem:pcextensions} that the domains be pairwise disjoint implies the following: we can distinguish $D_i$ as those elements in $\St(a_1)\setminus \St(\alpha_i)$ which are contained in $\St(\alpha_j)$ for every $j\neq i$.  The following corollary is immediate from this description.

\cor
Suppose $\Lambda$ contains $k+1$ vertices $a_1, \alpha_1, \dotsc, \alpha_k$ satisfying the following properties:
\enum
\item $a_1, \alpha_1, \dotsc, \alpha_k$ forms a clique.
\item $\St(a_1) \cup \bigcup \St(\alpha_i) = \Lambda$.
\item The sets $D_i = \left(\St(a_1) \setminus \St(\alpha_i)\right) \cap \bigcap_{j\neq i} \St(\alpha_j)$ are all nonempty.
\eenum
Define $\Gamma$ to be the graph obtained from $\Lambda$ by removing the vertices $\alpha_1, \dotsc, \alpha_k$ and any edge from $a_1$ to any $D_i$.  Then $W_{\Lambda}$ can be realized as the semi-direct product $W_{\Gamma} \rtimes H$, where $H \leq \Out^0(W_{\Gamma})$ is generated by the partial conjugations with acting letter $a_1$ and domains $D_i$.
\ecor

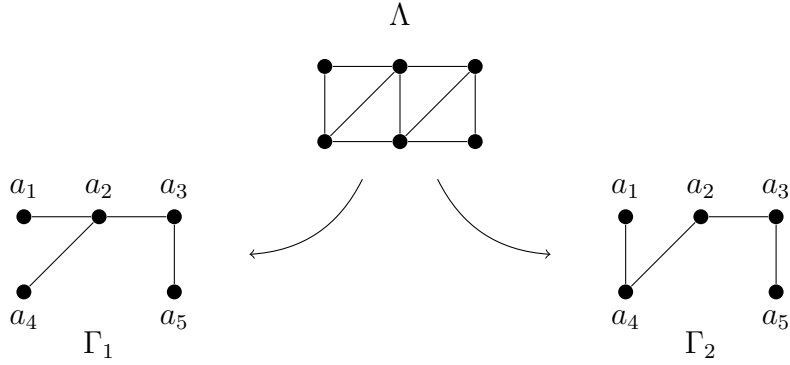
\begin{figure}[h]
\begin{tikzpicture}[scale=1.0]

\node (b1) [circle,fill,inner sep=2pt] at (-1,1) {};
\node (b2) [circle,fill,inner sep=2pt] at (0,1) {};
\node (b3) [circle,fill,inner sep=2pt] at (1,1) {};
\node (b4) [circle,fill,inner sep=2pt] at (-1,0) {};
\node (b5) [circle,fill,inner sep=2pt] at (0,0) {};
\node (b6) [circle,fill,inner sep=2pt] at (1,0) {};
\node at (0,1.7) {$\Lambda$};

\draw (b4) to (b1) to (b2) to (b3) to (b6) to (b5) to (b4) to (b2) to (b5) to (b3);

\begin{scope}[xshift=4cm,yshift=-2cm]
\node (b1) [circle,fill,inner sep=2pt,label=90:$a_1$] at (-1,1) {};
\node (b2) [circle,fill,inner sep=2pt,label=90:$a_2$] at (0,1) {};
\node (b3) [circle,fill,inner sep=2pt,label=90:$a_3$] at (1,1) {};
\node (b4) [circle,fill,inner sep=2pt,label=-90:$a_4$] at (-1,0) {};
\node (b6) [circle,fill,inner sep=2pt,label=-90:$a_5$] at (1,0) {};
\node at (0,-0.7) {$\Gamma_2$};

\draw (b1) to (b4) to (b2) to (b3) to (b6);
\end{scope}

\begin{scope}[xshift=-4cm,yshift=-2cm]
\node (b1) [circle,fill,inner sep=2pt,label=90:$a_1$] at (-1,1) {};
\node (b2) [circle,fill,inner sep=2pt,label=90:$a_2$] at (0,1) {};
\node (b3) [circle,fill,inner sep=2pt,label=90:$a_3$] at (1,1) {};
\node (b4) [circle,fill,inner sep=2pt,label=-90:$a_4$] at (-1,0) {};
\node (b6) [circle,fill,inner sep=2pt,label=-90:$a_5$] at (1,0) {};
\node at (0,-0.7) {$\Gamma_1$};

\draw (b1) to (b2) to (b3) to (b6);
\draw (b2) to (b4);
\end{scope}

\draw [->,bend left] (-0.5,-0.5) to (-2,-1.5);
\draw [->,bend right] (0.5,-0.5) to (2,-1.5);

\end{tikzpicture}
\caption{$W_{\Lambda} = W_{\Gamma_1} \rtimes \langle x \rangle = W_{\Gamma_2} \rtimes \langle y \rangle$, where $x, y$ act like the partial conjugations $x = \chi_{4, \{1\}}$ and $y = \chi_{2, \{1\}}$ on $\Gamma_1, \Gamma_2$, respectively.}
\label{fig:multipledecompositions}
\end{figure}

Theorem \ref{thm:mainextensiontheorem} yields an analogous corollary, since in each case they tell how to build the defining graph of the extension from the original defining graph, and the process is always reversible.  It is not uniquely reversible.  A given right-angled Coxeter group will, in general, have many semi-direct product decompositions.  As an example, consider the decompositions shown in Figure \ref{fig:multipledecompositions}.

\section{Details}
\label{sec:details}

In this section we explore the properties of the clique graph, the star poset, and the involution graph introduced in Section \ref{sec:summary}.  We present detailed proofs of these properties, including proofs establishing claims made in that section and the correctness of our collapsing algorithms.

\subsection{The Clique Graph and the Star Poset}
\label{sec:cliquegraph}

Recall that, given a graph $\Gamma$, we write $\Gamma_I$ for the intersections of maximal cliques in $\Gamma$.  We begin by establishing a correspondence between the \emph{maximal clique structure} of a graph $\Gamma$ and its clique graph $\Gamma_K$.  By \emph{maximal clique structure}, we mean that there is a bijection between the maximal cliques of $\Gamma$ and those of $\Gamma_K$ which respects intersections.

\prop
\label{prop:cliquegraphstructure}
Suppose $\Gamma$ is a finite graph with maximal cliques $\Gamma_1, \dotsc, \Gamma_r$.  For any subset $I \subseteq \{1, 2, \dotsc, r\}$, write 
\[
\Gamma_{I} = \bigcap_{i \in I} \Gamma_i.
\]
Similarly, write $\Gamma_{K, 1}, \dotsc, \Gamma_{K, s}$ for the maximal cliques of $\Gamma_K$, and write $\Gamma_{K, I}$ for the intersections of maximal cliques.  Then, possibly after reindexing:
\enum
\item $r = s$,
\item each $\Gamma_{K, J}$ contains at least one $J$-minimal vertex (namely, the vertex labeled by the clique $\Gamma_J$),
\item $\Gamma_{K, i} = (\Gamma_i)_K$, (that is, $(\Gamma_i)_K$ naturally injects as a labeled graph into $\Gamma_K$, and the image is precisely $\Gamma_{K, i}$),
\item $\Gamma_{K, I} =  (\Gamma_I)_K$, and
\item if $\Gamma_I$ is a clique of size $k$, then $\Gamma_{K, I}$ is a clique of size $2^k - 1$.
\eenum
\eprop
\pf
(1) For each maximal clique $\Gamma_i$ in $\Gamma$, there is a corresponding vertex $v_i$ in $\Gamma_K$. This vertex is adjacent only to vertices representing subsets of $\Gamma_i$ since $\Gamma_i$ is maximal, and so $v_i$ is contained in the unique maximal clique $\St(v_i)$ in $\Gamma_K$.  In particular, since each no $v_i, v_j$ can be in the same maximal clique of $\Gamma_K$, we have $r \leq s$. 




Conversely, each vertex of the maximal clique $\Gamma_{K, i}$ is labeled by some clique of vertices in $\Gamma$.  Since $\Gamma_{K, i}$ forms a clique, the collection of all vertices of $\Gamma$ which appear in the labels of vertices of $\Gamma_{K, i}$ must form a clique $\Lambda$ in $\Gamma$.  It is clear that $\Lambda$ is maximal, since $\Gamma_{K, i}$ is.  Thus $\Lambda = \Gamma_j$ for some $j$.  That is, $s\leq r$, establishing (1).  The description we have just given of the cliques in $\Gamma_K$ also establishes the correspondence in (3), and therefore in (4).


As noted in the claim, the clique $\Gamma_J$ forms a vertex of $\Gamma_K$.  It is straightforward to see that this vertex in $J$-minimal in $\Gamma_{K, J}$, establishing (2).



Finally, if $\Gamma_I$ is a clique of size $k$, then every non-empty subset of vertices induces a clique, and so corresponds to a vertex in $\Gamma_{K,I}$. There are $2^{k} - 1$ of these subsets, which correspond to $2^{k}-1$ vertices in $\Gamma_{K,I}$.
\epf

Let $\Gamma$ be a finite graph with maximal cliques $\Gamma_1, \dotsc, \Gamma_r$.  As before, write $\Gamma_I$ for the intersections of the maximal cliques, and suppose $|\Gamma_I| = k_I$.  Then
\begin{equation}
\label{eq:cliqueintersection}
\sum_{I \supsetneq J} (-1)^{|I\setminus J| + 1}k_I \leq k_{J}.
\end{equation}
This is a direct application of the inclusion-exclusion principle, since the left hand side of the inequality counts the number of vertices in $\Gamma_J \cap \bigcup_{i \notin J} \Gamma_i$ (while the right-hand side is, by definition, the total number of vertices in $\Gamma_J$).  We have therefore established that any clique graph must satisfy the Maximal Clique, Minimal Vertex, and Inclusion-Exclusion Conditions.  This gives one direction of the characterization theorem:

\cliquecharacterization*

If we are faced with some graph which we do not know to be a clique graph, we can check directly that the intersections of maximal cliques have sizes of the form $n_I = 2^{k_I} - 1$, and we can check directly that the system of integers $k_I$ satisfies the inclusion-exclusion inequalities.  Thus, determining whether a graph arises as a clique graph is reduced to checking a system of integer inequalities (once we establish the other direction of the theorem).



\begin{figure}[h]
\begin{tikzpicture}[scale=1.0]

\node (a) [circle,fill,inner sep=1pt,label=180:$a$] at (90:1) {};
\node (b) [circle,fill,inner sep=1pt,label=-45:$b$] at (210:1) {};
\node (c) [circle,fill,inner sep=1pt,label=45:$c$] at (330:1) {};

\node (d) [circle,fill,inner sep=1pt] at ($(120:1) + (90:1)$) {};
\node (e) [circle,fill,inner sep=1pt] at ($(60:1) + (90:1)$) {};

\node (f) [circle,fill,inner sep=1pt] at ($(-120:1) + (210:1)$) {};
\node (g) [circle,fill,inner sep=1pt] at ($(180:1) + (210:1)$) {};

\node (h) [circle,fill,inner sep=1pt] at ($(-60:1) + (330:1)$) {};
\node (i) [circle,fill,inner sep=1pt] at ($(0:1) + (330:1)$) {};

\draw (a) to (b) to (c) to (a);
\draw (a) to (d) to (e) to (a);
\draw (b) to (f) to (g) to (b);
\draw (c) to (h) to (i) to (c);

\end{tikzpicture}
\caption{The triangle $\{a, b, c\}$ forms a maximal clique which fails the Inclusion-Exclusion Condition.  This was essentially the feature of Example \ref{ex:nonracg} which prevented the group in that example from being right-angled Coxeter.}
\label{fig:triangleoftriangles}
\end{figure}
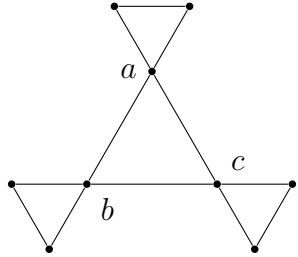

\ex
\label{ex:triangleoftriangles}
Consider the graph in Figure \ref{fig:triangleoftriangles}.  In this graph, all intersections of maximal cliques have sizes of the form $2^k - 1$, but the Inclusion-Exclusion Condition fails.  So the graph cannot arise as a clique graph.
\eex

We will establish the converse of Theorem \ref{thm:cliquecharacterization} by proving that, for any graph which satisfies the Maximal Clique, Minimal Vertex, and Inclusion-Exclusion Conditions, the proposed collapsing procedure of Theorem \ref{thm:collapsing} produces the desired output.  In order to evaluate the collapsing procedure, we must explore some properties of the star poset $\mathcal{P}(\Gamma)$.

\lem
Let $[v] \in \mathcal{P}(\Gamma)$.  Then the vertices
\[
S = \bigcup_{[v]\leq [w]} [w]
\]
form a clique in $\Gamma$.  If this clique is maximal, then $[v]$ is minimal in $\mathcal{P}(\Gamma)$.
\elem
\pf
If $w, w'\in S$ are any vertices, then $w \in \St(v) \subseteq \St(w')$, so $w$ and $w'$ are adjacent.  Thus $S$ forms a clique.

We now suppose $[v]$ is not minimal.  Then there is some $[w] < [v]$.  In particular, $w \notin S$, but $w\in \St(w) \subseteq \St(s)$ for any $s\in S$, hence $w$ is a vertex outside of $S$ adjacent to all of $S$.  Thus $S$ is not maximal.
\epf

\defi
For $[v] \in \mathcal{P}(\Gamma)$, we call the clique $S$ defined in the lemma the \emph{clique above $[v]$}.  We will use the notation $S_v$ if we need to keep track of the vertex $v$.
\edefi

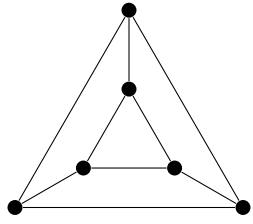
\begin{figure}[ht]
\begin{tikzpicture}[scale=0.7]

\node (v1) [circle,fill,inner sep=2pt] at (90:1) {};
\node (v2) [circle,fill,inner sep=2pt] at (210:1) {};
\node (v3) [circle,fill,inner sep=2pt] at (330:1) {};

\node (w1) [circle,fill,inner sep=2pt] at (90:2.5) {};
\node (w2) [circle,fill,inner sep=2pt] at (210:2.5) {};
\node (w3) [circle,fill,inner sep=2pt] at (330:2.5) {};

\draw (v1) to (v2) to (v3) to (v1);
\draw (w1) to (w2) to (w3) to (w1);
\draw (v1) to (w1);
\draw (v2) to (w2);
\draw (v3) to (w3);

\end{tikzpicture}
\caption{It is easy to check that each vertex is its own star equivalence class, and that these equivalence classes are pairwise not comparable.  In particular, each $[v]$ is minimal, and each $S_v = \{v\}$ is not a maximal clique.}
\label{fig:nonmaximalclique}
\end{figure}

The converse of the lemma (i.e., that minimality of $[v]$ implies maximality of $S_v$) is false in general.  A simple example is given in Figure \ref{fig:nonmaximalclique}.  However, we claim that the converse does hold for those $\Gamma$ which are clique graphs.  Namely:

\prop
\label{prop:minimal}
Suppose $\Gamma$ satisfies the Minimal Vertex Condition.  Then $[v]$ is a minimal element of $\mathcal{P}(\Gamma)$ if and only if $v$ is a minimal vertex of $\Gamma$.  In this case, $S_v$ is the unique maximal clique containing $v$.
\eprop
\pf
Suppose $v$ is a minimal vertex of $\Gamma$.  Then $\St(v)$ is the unique maximal clique containing $v$.  Since $S_v$ is a clique containing $v$, it is clear that $S_v \subseteq \St(v)$.  Conversely, if $x\in \St(v)$, then $\St(v) \subseteq \St(x)$, hence $[v]\leq [x]$ and $x\in S_v$.  Thus $\St(v) = S_v$ is maximal.  By the previous lemma, since $S_v$ is maximal, $[v]$ is minimal.

Conversely, suppose $v$ is not minimal.  Then $v$ is contained in the intersection of two distinct maximal cliques, $\Gamma_1$ and $\Gamma_2$.  Since $\Gamma_i$ are maximal cliques, they contain minimal vertices $w_i$.  By the above argument, $[w_i] \leq [v]$, and this must be a strict inequality since, e.g., $w_2 \in \St(v) \setminus \St(w_1)$.  Thus $[v]$ is not minimal.
\epf

\prop
\label{prop:intersections}
For any finite graph $\Gamma$ and $[v] \in \mathcal{P}(\Gamma)$, $S_v$ is an intersection of maximal cliques.
\eprop
\pf
Let $\Gamma_1, \dotsc, \Gamma_k$ be all the maximal cliques of $\Gamma$ containing $S_v$.  It is clear that $S_v \subseteq \bigcap \Gamma_i$.

Conversely, let $v'\in \bigcap \Gamma_i$ and suppose $v'\notin S_v$.  Since $\St(v) \nsubseteq \St(v')$, there is some $x\in \St(v)$ which is not in $\St(v')$.  In particular, since $\bigcup \Gamma_i \subseteq \St(v')$, we must have $x\notin \Gamma_i$ for any $i$.  By construction of $S_v$, we must have $x\in \St(w)$ for each $w\in S_v$.  Now $S_v \bigcup \{x\}$ forms a clique which contains $S_v$ and is not equal to $\Gamma_i$ for any $i$, contradicting our assumption that the list of $\Gamma_i$ contained all maximal cliques containing $S_v$.  So there can exist no such $v'$, hence $S_v = \bigcap \Gamma_i$, proving the claim.
\epf

We observe that the previous two propositions say the following in the case of clique graphs (which must satisfy the Minimal Vertex Condition):

\cor
\label{cor:minimal}
Suppose $\Gamma_K$ is a clique graph.
\enum
\item $[v]$ is minimal in $\mathcal{P}(\Gamma_K)$ if and only if $v$ is minimal in $\Gamma_K$, and
\item if $[v]$ is non-minimal, then $S_v$ is the intersection of maximal cliques (and therefore has size of the form $2^k - 1$).  In this case,
\[
S_v = \bigcap_{\substack{[w] \textrm{ minimal} \\ [w] \leq [v]}} S_w.
\]
\eenum
\ecor

This shows that the star poset also records information about the intersections of maximal cliques: any clique above $[v]$ is such an intersection.  Finally, we prove the converse.

\prop
\label{prop:yesallintersections}
Suppose $\Gamma_K$ is a clique graph.  Then any intersection of maximal cliques is equal to $S_v$ for some $v$.
\eprop
\pf
Since $\Gamma_K$ is a clique graph, it satisfies the Minimal Vertex Condition.  Let $\Gamma_{K, J}$ be any intersection of maximal cliques, and let $v\in \Gamma_{K, J}$ be a $J$-minimal vertex.  Without loss of generality, let $J$ be the maximal index set without changing the intersection.  In particular, $J$ is precisely the index set of all maximal cliques containing $v$, so that $\St(v) = \bigcup_{j\in J} \Gamma_{K,j}$.

We claim that $S_v = \Gamma_{K, J}$.  Let $u\in S_v$.  By definition of $S_v$, $[v] \leq [u]$, so
\[
\St(u) \supset \St(v) = \bigcup_{j\in J} \Gamma_{K,j}.
\]
That is, $u$ is adjacent to every vertex in $\Gamma_{K,j}$, for each $j\in J$.  Since each $\Gamma_{K,j}$ is a maximal clique, this shows $u\in \Gamma_{K,j}$ for each $j\in J$.  That is, $u\in \Gamma_{K, J}$.

Conversely, let $w\in \Gamma_{K, J}$.  Then $w$ is adjacent to all vertices in $\Gamma_{K,j}$ for $j\in J$, thus $\bigcup_{j\in J} \Gamma_{K,j} \subseteq \St(w)$.  That is, $\St(v) \subseteq \St(w)$, so $[v] \leq [w]$.  By definition of $S_v$, $w\in S_v$.
\epf

We note that the previous proof gives a nice description of the elements of each star equivalence class:

\cor
\label{cor:Jminimalclass}
Suppose $\Gamma_K$ is a clique graph.  Then any $[v] \in \mathcal{P}(\Gamma_K)$ consists precisely of the $J$-minimal vertices of $\Gamma_K$, where $J$ is the largest index set such that $v\in \Gamma_J$.
\ecor
\pf
Clearly, all $J$-minimal vertices for the same index set $J$ must have the same star (namely, $\bigcup_{j\in J} \Gamma_{K, j}$).  Conversely, suppose $v$ is $J$-minimal and $[v] = [w]$ for some $w$.  Then $w\in \Gamma_j$ for each $j\in J$, and $w\notin \Gamma_i$ for any $i\notin J$.  (Otherwise, all of $\Gamma_i$ would be in $\St(w)$, which we have assumed to be equal to $\St(v)$, a contradiction.)  Therefore, $w$ is $J$-minimal.
\epf

These results establish that, for a clique graph, the cliques above vertices are precisely the intersections of maximal cliques, and every intersection of maximal cliques is the clique above some vertex.  (This is not, in general, a bijective correspondence.  As remarked earlier, it may be that $\Gamma_{K, J} = \Gamma_{K, J'} = S_v$, where $J\neq J'$.)  In our collapsing algorithm to recover $\Gamma$ from $\Gamma_K$, we begin at the top of the poset (this is the deepest intersections of maximal cliques) and works downwards.  The previous proposition ensures that the algorithm examines every intersection of maximal cliques as it traverses every element in the poset structure.

We now with to prove the correctness of our collapsing procedure, which also establishes the other direction of Theorem \ref{thm:cliquecharacterization}.  Recall the procedure:

\collapsing*

We first must address a subtlety, namely, that we can carry out the choice in step 2 of the algorithm.

\prop
In step 2 of the collapsing procedure, $0 \leq k-k' \leq |[w]|$.  So we are able to choose an appropriate number of vertices from $[w]$ to add to $V$.
\eprop
\pf
The clique $S_w$ is some intersection of maximal cliques $\Gamma'_J$ by Corollary \ref{prop:minimal}.  From this clique, we have already chosen $k'$ vertices, and every vertex among those already chosen comes from a larger poset element, which is therefore a strictly smaller intersection of maximal cliques.  By the Inclusion-Exclusion Condition, the number of elements we could have chosen is at most $k_J = k$, hence $k'\leq k$.

Now $S_w = \left(\bigcup_{[w] < [v]} S_v\right) \cup [w]$.  Because 
\[
|S_w| = 2^k -1, \text{ and } \left|\bigcup_{[w] < [v]} S_v\right| \leq 2^{k'}-1, \text{ and } [w] \leq 2^{|[w]|},
\]
we have that $2^k -1 \leq 2^{k'} -1 + 2^{|[w]|}$.  Therefore $2^k \leq 2^{k'} + 2^{|[w]|}$. But $2^x + 2^y \leq 2^{x+y}$ for all pairs of positive integers $x, y$.  Thus $2^k \leq 2^{k'+|[w]|}$, and $k \leq k'+|[w]|$.
\epf

We also see that step 2 does not tell us explicitly which vertices of $[w]$ to add to $V$.  We claim this choice does not matter:

\prop
Given $\Gamma'$, if the procedure above does not return false, then the isomorphism type of the graph $\Gamma$ does not depend on the choices made in step 2 of the collapsing procedure.
\eprop
\pf
Without loss of generality, we will suppose our choices differ by a single vertex.  Suppose we are about to consider $[v]$ and have constructed the set $V$ thus far.  Let $v_1, \dotsc, v_{k+1} \in [v]$, where $k > 0$ is the number of vertices from $[v]$ which we must add to $V$.  Let
\begin{align*}
V_1 &= V \cup \{v_1, \dotsc, v_k\} \\
V_2 &= V \cup \{v_1,\dotsc, v_{k-1}, v_{k+1}\}.
\end{align*}
We observe that we can make all future choices the same (since we haven't changed the number of vertices we must pick from $[w]$ for any $[w]\leq [v]$), so that we create two final graphs $\Gamma_1$ and $\Gamma_2$ whose vertex sets differ only by switching $v_k$ and $v_{k+1}$.

We now claim that the resulting graphs $\Gamma_1$ and $\Gamma_2$ are isomorphic.  By the previous observation, the vertex sets of $\Gamma_1$ and $\Gamma_2$ differ only by switching $v_k$ and $v_{k+1}$.  So we can define a map $\varphi \colon \Gamma_1 \rightarrow \Gamma_2$ which sends each vertex other than $v_k$ to itself, and which sends $v_k$ to $v_{k+1}$.  We claim that $\varphi$ defines a graph isomorphism.  Clearly any adjacency relation not involving $v_k$ is preserved under $\varphi$.  Suppose $w$ is a vertex of $\Gamma_1$ adjacent to $v_k$.  Then $w\in \St(v_k) = \St(v_{k+1})$, so $w$ is adjacent to $v_{k+1}$.  Thus $\varphi$ is a graph homomorphism.  By the same argument, the analogous map $\psi \colon \Gamma_2 \rightarrow \Gamma_1$ is also a graph homomorphism, and the two maps are clearly inverses.  Hence $\Gamma_1$ is isomorphic to $\Gamma_2$.  The full result follows by induction.
\epf

This shows that the isomorphism type of an output graph $C(\Gamma')$ is determined.  However, a priori it could be the case that there are two graphs $\Gamma, \Lambda$ so that $\Gamma_K$ and $\Lambda_K$ are isomorphic, but the collapsing procedure applied to $\Gamma_K$ always outputs the isomorphism type $\Gamma$.  The following proposition says that the maximal clique structure of $\Gamma'$ determines the maximal clique structure of the output $C(\Gamma')$.  The theorem following the proposition establishes that the maximal clique structure (including information about the sizes of all intersections of maximal cliques) determines a graph up to isomorphism.  By Proposition \ref{prop:cliquegraphstructure}, any graph whose clique graph is $\Gamma'$ will have the same clique graph structure, and will therefore be isomorphic.  These results together show that the collapsing procedure outputs the unique graph $\Gamma$ up to isomorphism so that $\Gamma_K = \Gamma'$.

\prop
Let $\Gamma'$ be a finite graph satisfying the Maximal Clique, Minimal Vertex, and Inclusion-Exclusion Conditions.  In particular, this implies there is a system of integers $k_I$ so that $|\Gamma_I'| = 2^{k_I} - 1$.  Let $C(\Gamma') = \Gamma$.  Then the maximal cliques of $\Gamma$ correspond to the maximal cliques of $\Gamma'$, and $|\Gamma_I| = k_I$ for all $I$.
\eprop
\pf
By assumption, each $\Gamma_I'$ contains an $I$-minimal vertex $v_I'$.  We have $|\Gamma_I'| = 2^{k_I} - 1$, and the algorithm chooses exactly $k_I$ vertices from $S_{v_I'}$.  Corollary \ref{cor:minimal} implies that the maximal cliques in $\Gamma$ have sizes of the form $k_i$, and Proposition \ref{prop:yesallintersections} ensures that we have $|\Gamma_I| = k_I$ for all intersections of maximal cliques (since all intersections $\Gamma_I'$ occur as the clique above some element in the poset).
\epf

We have shown now that, if the algorithm returns any graph, then it returns a graph with a certain number of maximal cliques, and the intersections of the maximal cliques have certain sizes.  We now establish that a finite graph is determined up to isomorphism by the sizes of the intersections of maximal cliques.

\thm
Let $\Gamma, \Lambda$ be finite graphs.  Suppose both graphs have $r$ maximal cliques which may be indexed in such a way that, for all index sets $I\subset \{1, 2, \dotsc, r\}$, $|\Gamma_I| = |\Lambda_I|$.  That is, all intersections of maximal cliques have the same sizes in each graph.  Then there is an isomorphism $\varphi \colon \Gamma \rightarrow \Lambda$ which maps $\Gamma_i$ to $\Lambda_i$ for each $i$.
\ethm
\pf
We first claim that the poset structures $\mathcal{P}(\Gamma)$ and $\mathcal{P}(\Lambda)$ are the same, and the corresponding equivalence classes have the same sizes.  For each $v\in \Gamma$, let $J_v$ be the maximal index set so that $v\in \Gamma_{J_v}$.  Then $\St(v) = \bigcup_{j\in J_v} \Gamma_{j}$.  The equivalence class of $v$ consists of the $J_v$-minimal vertices of $\Gamma$ by Corollary \ref{cor:Jminimalclass}.  By assumption, $|\Gamma_{J_v}| = |\Lambda_{J_v}|$.  Moreover, the number of vertices which are in some further intersection is given by the inclusion-exclusion formula:
\[
\sum_{J \supsetneq J_v} (-1)^{|J\setminus J_v| + 1}|\Gamma_J| = \sum_{J \supsetneq J_v} (-1)^{|J\setminus J_v| + 1}|\Lambda_J|
\]
That is, the number of $J_v$-minimal vertices in $\Gamma$ and in $\Lambda$ is the same.  Since this is for any $v$, the sizes of star-equivalence classes of vertices in $\Gamma$ and $\Lambda$ are equal for every class.  Each equivalence class is represented by some index set $J$ (although not every index set represents a class).

An equivalence class represented by $J$ is smaller in the poset structure than another represented by $J'$ if and only if $J\subseteq J'$.  Since this holds in both $\Gamma$ and $\Lambda$, it follows that the poset structures are equivalent.

Now, we build a map $\varphi \colon \Gamma \to \Lambda$ by piecing together (arbitrary) bijections between each pair of corresponding equivalence classes.  We observe that, by construction, $\varphi([v]) = [\varphi(v)]$.

We also observe that $\Gamma_i$ is mapped to $\Lambda_i$ for each $i$: let $v\in \Gamma_i$, so that $i\in J_v$.  By construction, $\varphi(v)\in \Lambda_{J_v}$, which is an intersection of maximal cliques including $\Lambda_i$.  That is, $\varphi(v)\in \Lambda_i$.  It follows that $\varphi$ maps $\Gamma_I$ to $\Lambda_I$ for each $I$.

We must show that $\varphi$ preserves adjacency.  Suppose $v, w\in \Gamma$ are adjacent.  Then the edge $\{v, w\}$ extends to some maximal clique $\Gamma_i$.  Now $\varphi$ maps $\Gamma_i$ to $\Lambda_i$, so $\varphi(v)$ and $\varphi(w)$ are still adjacent.
\epf

This completes the proofs of Theorem \ref{thm:collapsing} and Theorem \ref{thm:cliquecharacterization}.

\subsection{Calculations in the Abelianization}
\label{sec:involutiongraph}

We now discuss the modifications to the collapsing procedure to make use of algebraic information.  Recall from the discussion in Section \ref{sec:summary} that, given a group $G$, we first form the involution graph $\Inv[G]$ and try to find a full system of representatives (i.e., a labeling of the vertices of $\Inv[G]$ which exhibit all commuting relations simultaneously).  If $\Inv[G]$ is a clique graph, the collapsing procedure will give a graph $\Gamma = C(\Inv[G])$ such that $\Gamma_K = \Inv[G]$.  Moreover, $\Gamma$ will carry the labels of the vertices chosen during the collapsing, so that the choice of which vertices to keep and which to omit is essentially the choice of which elements of $G$ will be the generators in a (hypothetical) right-angled Coxeter presentation.  For this reason, we must take care when choosing our generator vertices to avoid choosing group elements which have a nontrivial product relation.  We will now demonstrate a method of passing to the abelianization $\Ab{G}$ to determine product relations using straightforward calculations.

Suppose we are given a finitely presented group
\[
G = \langle s_1, \dotsc, s_m \mid r_1, \dotsc, r_k\rangle.
\]
Recall that, for $g\in G$, we will write $\overline{g}$ for the image of $g$ in the abelianization.  A presentation for $\Ab{G}$ is given by
\[
\Ab{G} \cong \left\langle \overline{s_1}, \overline{s_2}, \dotsc \overline{s_m} \mid \overline{r_1}, \overline{r_2}, \dotsc \overline{r_k}, [\overline{s_i}, \overline{s_j}] \textrm{ for } 1 \leq i, j \leq m \right\rangle.
\]
Writing the group operation additively in $\Ab{G}$, we can write the relations as linear combinations of the generators with integer coefficients:
\[
\overline{r_i} = a_{i,1} \overline{s_1} + a_{i,2} \overline{s_2} + \dotsb + a_{i,m} \overline{s_m}
\]
We collect the coefficients $(a_{i,j})$ into a $k \times m$ matrix $R$, called the \emph{relations matrix for $\Ab{G}$}.

We briefly recall the \emph{Smith normal form}: given the integer $k\times m$ matrix $R$, there exist $k\times k$ and $m\times m$ invertible matrices $P, Q$ and a diagonal matrix $S$ so that $R = PSQ$, and the diagonal elements of $S$ are $\alpha_1, \dotsc, \alpha_r, 0, \dotsc, 0$ such that $\alpha_i \mid \alpha_{i+1}$.  $S$ is called the \emph{Smith normal form of $R$}.

Interpreting $S$ as the relation matrix for a presentation, we have that $\Ab{G}$ is in a canonical form as a direct product of cyclic groups. Normal forms are immediate and computations in $\Ab{G}$ are much easier. Moreover, we now have an effective quotient map from $G \to \Ab{G}$ in this canonical form. Namely, for any $g \in G$ with $g = \prod s_j$, we have $\overline{g} = \sum \overline{s_j} = \sum_{i=1}^{m} b_i \overline{s_i}$. The vector-matrix product $\pmat b_1 & b_2 & \dotsb & b_m \epmat Q$ will give the coefficients of $\overline{g}$ in the Smith normal form presentation of $\Ab{G}$. This makes product relations easy to compute.

We now apply this method to show that, step 2 of the collapsing procedure, we can avoid nontrivial product relations.

\propchoice*

\pf
In step 2 of our collapsing procedure, we consider an equivalence class $[w]$ of $\Inv[W_{\Gamma}]$ and the clique above it, $S_w$, where $|S_w| = 2^k - 1$ for some $k$. If $(W_{\Gamma},S)$ is a right-angled Coxeter system for $W_{\Gamma}$ and the labels are distinct, pairwise commuting involutions, then $H = S_w \cup \{e\}$ is a finite subgroup isomorphic to $(\Z/2\Z)^k$:
the elements of $S_w$ are all involutions which pairwise commute.  Any product $g$ of these elements is an involution and commutes with all other elements of $S_w$ (so it is connected to all of $S_w$).  Moreover, any $h$ which commutes with all of $S_w$ commutes with any product of elements in $S_w$ (namely $g$), and so $g$ is contained in any maximal clique containing all of $S_w$. Since $S_w$ is an intersection of maximal cliques and $g$ is in all of these cliques, therefore $g\in S_w$.  So $H$ is a subgroup.



By Corollary \ref{cor:inj}, this subgroup projects injectively as a vector subspace into $\Ab{W_{\Gamma}}$. Inductively, we assume that there exists a choice of a right-angled Coxeter system $(W_{\Gamma},S)$ such that $\overline{V}$ is a set of \emph{standard basis elements} for $\Ab{(W_{\Gamma},S)} \cong (\Z / 2\Z)^k$, i.e., each element has only one non-zero component in the representation for the abelianization given by our choice of right-angled Coxeter system $(W_{\Gamma},S)$. (The base case is $\overline{V} = \emptyset$ and any choice of $(W_{\Gamma},S)$.) 

It follows that $\overline{V \cap S_w}$ is a linearly independent set in the $\mathbb{Z}/2\mathbb{Z}$-vector space $\Ab{W_{\Gamma}}$. We can then choose $k - k'$ labels in $S_w - V$ to extend this linearly independent set to a basis $\overline{B}$ of $\langle \overline{S_w} \rangle$. (It's possible that $k - k' = 0$.) Since $H$ projects injectively, choosing a basis for $\langle \overline{S_w} \rangle$ is the same as choosing a basis for $\langle S_w \rangle$. We need to show that $\overline{V} \cup \overline{B}$ is linearly independent as well. 

To clarify, we are now keeping track of \emph{two} different representations of the abelianization. $\Ab{W_{\Gamma}}$ is the form calculated from the Smith Normal Form in step 0 of the procedure, and $\Ab{(W_{\Gamma},S)}$ is the form wherein each element of $\overline{V}$ is a standard basis element. We will show that such a form must exist if $W_{\Gamma}$ is a right-angled Coxeter group, but it will not be directly computable during the procedure itself. The existence of this form will be used to show that $\emph{any}$ choice of $\overline{B}$ during our procedure will result in no non-trivial product relations. 

Since $H$ is a finite subgroup of $(W_{\Gamma},S)$, it is conjugate to a special subgroup: $g H g^{-1} = \langle a_1, a_2, \dotsc, a_k \rangle$ for $\{a_1, a_2, \dotsc, a_k\} \subseteq S$. Consider $b \in B \subseteq S_w$. Reordering the vertices of $S$ if necessary, then $g b g^{-1} = a_1 a_2 \dotsm a_m$ in $(W_{\Gamma},S)$. By the deletion condition of right-angled Coxeter groups (see, for example, \cite{Davis}), a product of distinct, commuting generators of $(W_{\Gamma},S)$, $c_1 c_2 \dotsm c_{\ell}$, commutes with $a_1 a_2 \dotsm a_m$ if and only if $c_j$ commutes with $a_i$ for each $i,j$. In particular, $[b] = [a_1 a_2 \dotsm a_m] \leq [a_i]$ for each $1 \leq i \leq m$.

Suppose that $[b] \lneq [a_i]$ for each $i$. Then the procedure has already considered $[a_i]$, and a subset of $\overline{V}$ is a basis for $\langle \overline{S_{a_i}} \rangle$, which contains $\overline{a_i}$. But by our inductive hypothesis, $\overline{V}$ is a set of standard basis elements relative to $\Ab{(W_{\Gamma},S)}$ and since $a_i \in S$, $\overline{a_i}$ is also a standard basis element. So the only way that $\overline{a_i} \in \langle \overline{V} \rangle$ is if $\overline{a_i} \in \overline{V}$ (and so by injectivity $g^{-1} a_i g \in V$). Thus, $b = g^{-1} a_1 a_2 \dotsm a_m g \in \langle V \rangle$ and so would \emph{not} be chosen by the procedure to linearly extend $\overline{V}$. 

Therefore, there must be some $i$ such that $[b] = [a_i]$. By reordering the vertices of $S$ if necessary, $[b] = [a_1]$. But then $g b g^{-1} = a_1 a_2 \dotsm a_m$ and $a_1$ are involutions that commute with exactly the same involutions, and so the following is an involutive automorphism (in fact a transvection) of $(W_{\Gamma},S)$: 
\begin{align*}
&\varphi: W_{\Gamma} \to W_{\Gamma} \\
\varphi(a_j) =
  &\begin{cases}
  a_1 a_2 \dotsm a_m  & \text{if } j=1 \\
  a_j       & \text{otherwise }
  \end{cases}
\end{align*}

Now, $(W_{\Gamma}, \varphi(S))$ is also a right-angled Coxeter system for $W_{\Gamma}$ with the exact same generators except for swapping $a_1$ and the product $a_1 a_2 \dotsm a_m$. $\overline{V}$ is a set of standard basis elements \emph{not} including $\overline{a_1}$ and so is unchanged under the induced map $\overline{\varphi}: \Ab{(W_{\Gamma},S)} \to \Ab{(W_{\Gamma},\varphi(S))}$. Alternatively, $\varphi(b) = \varphi(g)^{-1} a_1 \varphi(g)$ and so $\overline{\varphi(b)} = \overline{a_1}$. So if we let $(W_{\Gamma},S') = (W_{\Gamma}, \varphi(S))$ be our new right-angled Coxeter system and $V' = V \cup \{b\}$ be our new subset of labels from our chosen full set of representatives of $\Inv[W_{\Gamma}]$, then the inductive hypothesis is still satisfied. In particular, in our Smith Normal Form $\Ab{W_{\Gamma}}$, $\overline{V'}$ is still linearly independent. 

For each $b \in B$, we can perform this procedure in succession making sure that for each $b$, we choose different $a_i$ such that $[b] = [a_i]$. If at any point this were not possible, it would mean that there was some $b_n = g^{-1} a_1 a_2 \dotsm a_m g$ (in the updated system $(W_{\Gamma},S')$ with $V'$) such that each $a_j$ either satisfies:
\begin{enumerate}
\item $[b_n] \lneq [a_j]$ in which case $\overline{a_j} \in \overline{V'}$ from a previous step in the procedure, or
\item $\overline{b_l} = \overline{a_j}$ for some $l < n$ in which case $\overline{a_j} \in \overline{V'}$ from a previous element of the basis.
\end{enumerate}
In either case, since all of the $\overline{a_j} \in \overline{S_w}$, this would give a linear dependence in $\overline{S_w}$ among $\overline{B}$, which contradicts its choice as a basis. 

Thus, by induction on both elements of the poset, and then within each class on the elements of each chosen basis, it will always be the case that $\overline{V}$ will consist of elementary basis elements in $\Ab{(W_{\Gamma},S)}$ for some choice of system $(W_{\Gamma},S)$. Since every generator $a_i$ of $S$ is in $S_{a_i}$, $\overline{a_i} \in \langle \overline{V} \cap \overline{S_{a_i}} \rangle$, but since $\overline{V}$ are all elementary basis vectors, it must be that $\overline{a_i} \in \overline{V}$. Thus, at the end of the procedure, $\overline{V}$ will always be the full standard basis for some system $\Ab{(W_{\Gamma},S)}$, and in particular, $\overline{V}$ will always be a basis of $\Ab{W_{\Gamma}}$. 

Any non-trivial product relation among the elements of $V$ would induce a linear dependence among their images in $\Ab{W_{\Gamma}}$. But since $\overline{V}$ is a basis, this can never happen. 
\epf

Finally, we prove the proposition that allows us to hypothetically build edges in the involution graph of a given group by doing calculations in the abelianization:

\invedges*
\pf
Let $(W_{\Gamma},S)$ be a right-angled Coxeter system for $W_{\Gamma}$, and let $[x]$ and $[y]$ be conjugacy classes of involutions. Since $x$ and $y$ are involutions in a right-angled Coxeter group, they are each conjugate to a product of commuting generators.   So there exists $a_1, a_2, \dotsc, a_n \in S$,  $b_1, b_2, \dotsc b_m \in S$, $g, h \in W_{\Gamma}$ such that
\begin{align*}
gxg^{-1} &= a_1 a_2 \dotsm a_n \\
hyh^{-1} &= b_1 b_2 \dotsm b_m
\end{align*}
where all of the $a_i$ pairwise commute, and all of the $b_j$ pairwise commute. Consider the product
\[
w = a_1 a_2 \dotsm a_n b_1 b_2 \dotsm b_m = c_1 c_2 \dotsm c_k,
\]
where the $c_{\ell}$ are the generators that appear among either the $a_i$ or the $b_j$ but not both. (The ones that appear in both cancel with each other since they can be brought to the front or back of their respective words.) In the abelianization $\Ab{W_{\Gamma}}$, $\overline{x} = \overline{a_1 a_2 \dotsm a_n}$, $\overline{y} = \overline{b_1 b_2 \dotsm b_m}$, and $\overline{w} = \overline{c_1 c_2 \dotsm c_k}$. 

Now suppose that $[x]$ and $[y]$ are connected by an edge in $\Inv[W_{\Gamma}]$. That means that some conjugates of $x$ and $y$ commute. This implies that the product of those conjugates, $z$, is an involution. But then in $\Ab{W_{\Gamma}}$, $\overline{z} = \overline{x} \overline{y}$. 

Conversely, suppose that there exists an involution, $z$ such that $\overline{z} = \overline{x} \overline{y}$. Since $z$ is an involution, it must be conjugate to a product of distinct, commuting generators, each of which is mapped to its corresponding generator of $\Ab{W_{\Gamma}}$ and so can be recovered directly from $\overline{z}$. Thus, these generators must be exactly the $c_l$, and so they each pairwise commute. In particular, $w$ is an involution, and $gxg^{-1}$ and $hyh^{-1}$ commute. Thus, $[x]$ and $[y]$ should be connected by an edge in $\Inv[W_{\Gamma}]$. 

\epf

We have now established the correctness of our right-angled Coxeter recognition procedure:

\algebraiccollapsing*

\section{Further Research}
\label{sec:future}

While we have used our decision procedure to successfully establish both positive and negative identification of right-angled Coxeter presentations among extensions of right-angled Coxeter groups, much work remains to be done. One might hope to eventually characterize all subgroups $H \leq \Out^0(W_{\Gamma})$ (or $H\leq \Aut(W_{\Gamma})$) such that $W_{\Gamma} \rtimes H$ is right-angled Coxeter.  We note that subgroups $H \leq \Out^0(W_{\Gamma})$ are not necessarily generated by partial conjugations (they may be generated by products of partial conjugations).  Even if we only considered those $H$ generated by partial conjugations, we could not extend Lemma \ref{lem:pcextensions} by induction.  If $x, y$ are two commuting partial conjugations of $W_{\Gamma}$, then 
\[
W_{\Gamma} \rtimes \langle x, y \rangle \cong \left(W_{\Gamma} \rtimes \langle x \rangle\right) \rtimes \langle y \rangle,
\]
however, $y$ may not act on $W_{\Gamma} \rtimes \langle x \rangle$ as a partial conjugation (it will generally act as a product of partial conjugations).  Theorem \ref{thm:mainextensiontheorem} extends the lemma by induction somewhat, but we have many more examples of right-angled Coxeter extensions which are not covered by this theorem.  More work is required for a complete characterization.

As in Section \ref{sec:decompositionresults}, following a characterization of extensions $W_{\Gamma} \rtimes H$ which are right-angled Coxeter, we would also gain insight into semi-direct product decompositions of right-angled Coxeter groups.  Given a graph $\Lambda$, we could hope to obtain a complete list of graph features which identify $W_{\Lambda}$ as $W_{\Gamma} \rtimes H$, where $H\leq \Out^0(W_{\Gamma})$.  (We observe that this would not identify all semi-direct product decompositions of right-angled Coxeter groups.  There are certainly decompositions which are not of this form.)

We strongly suspect that, whenever $H \leq \Out^0(W_{\Gamma})$ is isomorphic to $D_{\infty}$, then $W_{\Gamma} \rtimes H$ is not right-angled Coxeter.  The first example of Section \ref{sec:negativeresults} is of this form.  Much of the argument in that example rests on using a normal form to establish that the given list of classes of involutions is complete.  A general proof would require substantially more work to prove that we can accurately build the involution graph in the general case. 

In particular, in the case of universal right-angled Coxeter groups (those whose defining graphs have no edges), the outer automorphism groups act on a contractible simplicial complex called McCullough--Miller space \cite{Piggott}. This space is analogous to Culler--Vogtmann Outer space for the case of free groups \cite{Culler}, and we can use the action to classify all conjugacy classes of involutions in the outer automorphism groups. An analogous structure does not currently exist for the outer automorphism group of a general right-angled Coxeter group, and such a theory would need to be developed in order to construct the involution graph and confirm our conjecture.

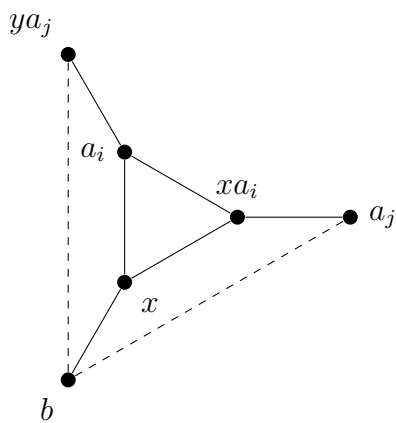
\begin{figure}[h]
\begin{tikzpicture}[scale=1.0]

\node (ai) [circle,fill,inner sep=2pt,label=180:$a_i$] at (120:1) {};
\node (x) [circle,fill,inner sep=2pt,label=-45:$x$] at (240:1) {};
\node (xai) [circle,fill,inner sep=2pt,label=90:$xa_i$] at (0:1) {};
\node (yaj) [circle,fill,inner sep=2pt,label=120:$ya_j$] at (120:2.5) {};
\node (b) [circle,fill,inner sep=2pt,label=240:$b$] at (240:2.5) {};
\node (aj) [circle,fill,inner sep=2pt,label=0:$a_j$] at (0:2.5) {};

\draw (ai) to (x) to (xai) to (ai);
\draw (ai) to (yaj);
\draw (x) to (b);
\draw (xai) to (aj);

\draw [dashed] (b) to (yaj);
\draw [dashed] (b) to (aj);

\end{tikzpicture}
\caption{The dotted lines represent edges that may be present in some cases.}
\label{fig:Dinftyproblem}
\end{figure}

Nevertheless, we can provide the following heuristic about what ought to go wrong in such an extension. Consider, for simplicity, a $D_{\infty}$ generated by two non-commuting partial conjugations.  If $x = \chi_{i, D}$ and $y = \chi_{j, E}$ are the partial conjugations, let $b$ be any vertex other than $a_j$ which is outside $\St(a_i) \cup D$.  Then Figure \ref{fig:Dinftyproblem} shows part of the involution graph of the extension.

In the figure, the edge from $b$ to $ya_j$ will be present if $b\in E$; the edge from $b$ to $a_j$ will be present if $b\in \St(a_j)$.  The figure as drawn so far cannot be a clique graph, because the central triangle is a maximal clique which does not satisfy the Inclusion-Exclusion Condition.  But even if other vertices were present which could turn the central triangle into a 7-clique (or larger) so that the condition would be satisfied, the collapsing procedure would need to choose all three vertices $x, a_i, xa_i$, which are not linearly independent in the abelianization.  However, this only establishes that the given pattern of labeling vertices in the involution graph---a pattern which has produced full systems of labels in all other examples so far---does not give an isomorphism to a right-angled Coxeter group in this case.  We have not sufficiently established that the extension could not have any right-angled Coxeter presentation.

\section*{Acknowledgements}

This work was partially supported by a grant from the Simons Foundation (\#317466 to Adam Piggott).  The authors would also like to thank Mauricio Gutierrez and the anonymous referee for a careful reading and helpful comments and suggestions.

\newpage
\bibliographystyle{alpha}
\bibliography{racg.bib}

\newpage
\todos

\end{document}